%% file: On_the_Correspondence_Between_Integer_Sequences_and_Vacillating_Tableaux.tex
\documentclass[11pt]{amsart} 
\usepackage[utf8]{inputenc}
\usepackage{fullpage}

\usepackage[english]{babel}
\usepackage[left=1 in, right=1 in, top=1.4 in, bottom=1.3 in]{geometry}
\usepackage{amsmath}
\usepackage{amssymb}
\usepackage{mathtools}
\usepackage[utf8]{inputenc}
\usepackage{color}
\usepackage{tikz}
\usetikzlibrary{patterns}
\usetikzlibrary{graphs,arrows.meta}
\usepackage[shortlabels]{enumitem}
\usepackage{graphicx}
\usepackage{caption}
\usepackage{subcaption}
\usepackage{hyperref}
\usepackage{tikz}
\usepackage{ytableau}
\usepackage{tabularx}
\usepackage{mathtools}
\usepackage[normalem]{ulem}

\theoremstyle{definition}
\newtheorem{definition}{Definition}
\newtheorem{proposition}[definition]{Proposition}
\newtheorem{lemma}[definition]{Lemma}
\newtheorem{theorem}[definition]{Theorem}
\newtheorem{corollary}[definition]{Corollary}
\newtheorem{example}[definition]{Example} 

\newcommand{\bsy}{\boldsymbol}

\newcommand{\vt}{vacillating tableau}
\newcommand{\vtx}{vacillating tableaux}

\usepackage{soul}

 \usepackage[colorinlistoftodos]{todonotes}


\title{
On the Correspondence Between Integer Sequences and Vacillating Tableaux}

\author{Zhanar Berikkyzy}
\author{Pamela E.~Harris}
\author{Anna Pun}  
\author{Catherine Yan}
\author{Chenchen Zhao} 

\address[Berikkyzy]{Department of Mathematics, Fairfield University, Fairfield, CT 06824 } 
\email{\textcolor{blue}{\href{mailto:berikkyzy@fairfield.edu}{berikkyzy@fairfield.edu}}}

\address[Harris]{Department of Mathematical Sciences, University of Wisconsin Milwaukee-Milwaukee, Milwaukee, WI 53211} 
\email{\textcolor{blue}{\href{mailto:peharris@uwm.edu}{peharris@uwm.edu}}}

\address[Pun]{Department of Mathematics, CUNY Baruch College, New York, NY 10010} 
\email{\textcolor{blue}{\href{mailto:anna.pun@baruch.cuny.edu}{anna.pun@baruch.cuny.edu}}}

\address[Yan]{Department of Mathematics, Texas A\&M University, College Station, TX 77843}
\email{\textcolor{blue}{\href{mailto:huafei-yan@tamu.edu}{huafei-yan@tamu.edu}}}

\address[Zhao]{Department of Mathematics, University of Southern California, Los Angeles, CA 90089} 
\email{\textcolor{blue}{\href{mailto:zhao109@usc.edu}{zhao109@usc.edu}}}

\date{\today}

\begin{document}

\begin{abstract}
A fundamental identity in the representation theory of the partition algebra is 
$n^k = \sum_{\lambda} f^\lambda m_k^\lambda$ for $n \geq 2k$, 
where $\lambda$ ranges over integer partitions of $n$, $f^\lambda$ is the number of standard Young tableaux of shape $\lambda$, and $m_k^\lambda$ is the number of vacillating tableaux of shape $\lambda$ and length $2k$. 
Using a combination of RSK insertion and jeu de taquin,  Halverson and Lewandowski constructed a bijection $DI_n^k$ that maps each integer sequence in $[n]^k$ to a pair of tableaux of the same shape,  where one is a standard Young tableau and the other is a vacillating tableau. In this paper, we study the fine properties of Halverson and Lewandowski's bijection
and explore the correspondence between integer sequences and the \vtx\ 
via the map $DI_n^k$ for general integers $n$ and $k$. 
In particular, we characterize the integer sequences $\bsy{i}$ whose corresponding shape, $\lambda$, in the image $DI_n^k(\bsy{i})$, satisfies $\lambda_1 = n$ or $\lambda_1 = n -k$.

\end{abstract}
\maketitle

\section{Introduction} 
A fundamental identity in the representation theory of the partition algebra is 
\begin{equation} \label{eq: Identity1} 
    n^k=\sum_{\lambda\vdash n} f^\lambda m_{k}^\lambda, 
\end{equation}
for $n \ge 2k$, where the sum is over integer partitions $\lambda$ of $n$, $f^\lambda$ is the number of standard Young tableaux (SYT) of shape $\lambda$, and $m_k^\lambda$ is the number of vacillating tableaux of shape $\lambda$ and length $2k$. 
This identity reflects a combinatorial analogue of the Schur-Weyl duality between the symmetric group algebra and the partition algebra.
Halverson and Lewandowski \cite{HL05} constructed an elegant  bijective proof  of  Identity \eqref{eq: Identity1}. 
Using a deletion-insertion algorithm based on the Robinson-Schensted-Knuth (RSK) insertion algorithm and jeu de taquin (jdt), Halverson and Lewandowski associate to each integer sequence $\bsy{i} \in [n]^k$, where $[n]\coloneqq\{ 1, 2, \ldots, n\}$,  a pair of tableaux of the same shape, one being a standard Young tableau and the other a vacillating tableau. 
The purpose of the present work is to study the properties of Halverson and Lewandowski's algorithm and to explore the correspondence between integer sequences and the shape of the associated tableaux. 

We begin by providing necessary definitions. 
A partition of a positive integer $n$ is a sequence  $\lambda=(\lambda_1, \lambda_2,\ldots, \lambda_t)$ of integers such that $\lambda_1 \geq \lambda_2 \geq \cdots \geq \lambda_t >0$ and $|\lambda|\coloneqq\lambda_1+ \cdots +\lambda_t =n$. 
We say that the size of $\lambda$ is $n$ and denote it as  $\lambda \vdash n$. 
In addition, we let the empty partition $\emptyset$ be the only integer partition of $0$.

A partition $\lambda$  can be 
visually represented by the Young diagram whose $j$-th row contains  $\lambda_j$ boxes.   
We adopt the English notation in which the diagrams are left justified with the first row at the top. 
A \emph{Young tableau} of shape $\lambda$ is an array obtained by filling each box of the Young diagram of $\lambda$ with an integer. 
A Young tableau is \emph{semistandard} if the entries are  weakly increasing in every row and strictly increasing in every column. 
A semistandard Young tableau (SSYT) is \emph{partial} if the entries are all distinct, in which case we call it a \emph{partial tableau}. 
 The \emph{content} of a Young tableau $T$, denoted as content$(T)$,  is the multiset of all the entries in $T$.
If the content of an SSYT $T$ with shape $\lambda \vdash n$ is exactly $[n]$, then we say $T$ is 
  a \emph{standard Young tableau} (SYT).  
  Throughout this paper, we use $f^\lambda$ to denote the number of SYTs of shape $\lambda$. 
  By convention, we set $f^\emptyset = 1$.

Below is the definition of a vacillating tableau as 
 introduced in \cite{HL05}. 
To emphasize the dependency on the parameter $n$, 
we call it an $n$-vacillating tableau. 

\begin{definition}[\cite{HL05}]
Let $n>0$ and $k\ge 0$ be integers.  An \emph{$n$-vacillating tableau of shape $\lambda$ and length $2k$} is a sequence of integer partitions 
\[
\Gamma=((n)=\lambda^{(0)}, \lambda^{(\frac{1}{2})}, \lambda^{(1)}, \lambda^{(1\frac{1}{2})}, \ldots, \lambda^{(k-\frac{1}{2})}, \lambda^{(k)}=\lambda) 
\]
such that for each integer $j=0, 1, \ldots, k-1$,  
\begin{enumerate}[(a)]
    \item $\lambda^{(j)} \supseteq  \lambda^{(j+\frac12)}$ and 
    $| \lambda^{(j)} / \lambda^{(j+\frac12)}|=1$, 
    \item $\lambda^{(j+\frac12)}\subseteq \lambda^{(j+1)}$ and 
    $|\lambda^{(j+1)}/ \lambda^{(j+\frac12)}|=1$. 
    \end{enumerate} 
Note that the above conditions imply that   $\lambda^{(j)} \vdash n$ and $\lambda^{(j+\frac12)} \vdash (n-1)$ for each integer $j$. 
\end{definition} 
    
Let $\mathcal{VT}_{n,k}(\lambda)$ denote the set of $n$-vacillating tableaux of shape $\lambda$ and  length $2k$. Following \cite{HL05}, we use 
$m_{k}^\lambda$ 
for the cardinality of $\mathcal{VT}_{n,k}(\lambda)$.
From the definition and using induction, we have that $\lambda^{(j)}$ has at least $n-j$ boxes in the first row. 
Hence, $m_{k}^\lambda$ is non-zero only when the first entry $\lambda_1$ of $\lambda$ 
satisfies $\lambda_1 \geq n-k$.

There is an equivalent notion of vacillating tableau that does not depend on the parameter $n$, which was independently introduced 
by Chen et.al. in \cite{CDDSY07}. 
To distinguish from the $n$-vacillating tableau, we call this second notion the \emph{simplified vacillating tableau}. 
Given an integer partition $\lambda=(\lambda_1, \lambda_2, \ldots, \lambda_t) \vdash n$, let 
$\lambda^*=(\lambda_2, \ldots, \lambda_t) \vdash (n-\lambda_1)$.
For an $n$-vacillating tableau $\Gamma= ( \lambda^{(j)}: j=0, \frac{1}{2}, 1, 1\frac{1}{2}, \ldots, k)$ in $\mathcal{VT}_{n,k}(\lambda)$, the \emph{simplified vacillating tableau $\Gamma^*$} is the sequence 
$( \mu^{(j)}: \mu^{(j)}= (\lambda^{(j)})^* \text{ for  } j=0, \frac{1}{2}, 1, 1\frac{1}{2}, \ldots, k)$. 
One can also define the simplified vacillating tableau in terms of integer partitions as in Definition 2.  

\begin{definition}[\cite{CDDSY07}]
A \emph{simplified vacillating tableau} $\Gamma^*$ of shape $\mu$  and length $2k$  is a sequence $(\mu^{(j)}: j = 0, \frac{1}{2}, 1, 1\frac{1}{2}, \ldots, k)$ of integer partitions such that $\mu^{(0)} = \emptyset$, $\mu^{(k)} = \mu$, and for each integer $j = 0, 1, \ldots, k-1$,
\begin{enumerate}[(a)]
\item $\mu^{(j)} \supseteq  \mu^{(j+\frac12)}$ and 
    $| \mu^{(j)} / \mu^{(j+\frac12)}|=0$ or $1$, 
    \item $\mu^{(j+\frac12)}\subseteq \mu^{(j+1)}$ and 
    $|\mu^{(j+1)}/ \mu^{(j+\frac12)}|=0$ or $1$. 
\end{enumerate}
 \end{definition} 
 
Denote by $\mathcal{SVT}_k(\mu)$ the set of simplified vacillating tableaux of shape $\mu$ and length $2k$. 
If $n \geq 2k$, it is clear that $\Gamma \leftrightarrow \Gamma^*$ is a bijection between $\mathcal{VT}_{n,k}(\lambda)$ and $\mathcal{SVT}_k(\lambda^*)$. 

The set of simplified vacillating tableaux has a close tie with discrete structures. 
In \cite{CDDSY07} a bijection is constructed between 
set partitions on $[k]$ and 
simplified vacillating tableaux of the empty shape and length $2k$.
This bijection reveals the symmetry between 
maximal crossings and maximal nestings of set partitions. 
 In general, a simplified vacillating tableau of shape $\lambda$ can be represented by a pair consisting of an arc diagram and a standard Young tableau
 of shape $\lambda$,  which yields many combinatorial formulas and identities. 
See, for example,  \cite{BHPYZ1, BHPYZ2, CDDSY07, DuYan24}.  
The same idea was further extended to fillings of Ferrers shapes, stack polyominoes, and moon polyominoes, under the language of growth diagrams, to describe the properties of maximal increasing and decreasing chains in those fillings. 
For more on those connections, see  \cite{GP2020,Kratten06, Kratten23, Rubey11}.

On the other hand, the combinatorial properties of $n$-vacillating tableaux for general $n$ and $k$ are not well-understood. 
As mentioned previously, the notion of $n$-\vt\ was introduced in \cite{HL05} to prove  Identity \eqref{eq: Identity1}.   
For any integer sequence $\bsy{i}=(i_1, i_2, \ldots, i_k) \in [n]^k$, Halverson and Lewandowski define an iterative delete-insert process, which we denote by $DI_{n}^k$, that  associates to $\bsy{i}$ a pair $DI_n^k(\bsy{i})=(P_{\lambda}(\bsy{i}), \Gamma_{\lambda}(\bsy{i}))$, where $\lambda$ is an integer partition of $n$, $P_{\lambda}(\bsy{i})$ is an SYT  of shape $\lambda$ and $\Gamma_{\lambda}(\bsy{i})$ is an $n$-\vt\ of shape $\lambda$. 
The exact definition of $DI_n^k(\bsy{i})$ is given in Section~\ref{Sec:algorithm}. 
The combinatorial properties of $DI_n^k$ have not been thoroughly studied. For example, 
 there is no known growth diagram description of $DI_n^k$ yet. 
 
The objective of this paper is to understand  the correspondence between integer sequences and the \vtx\ 
via the map $DI_n^k$ for general integers $n$ and $k$. 
In particular, we are interested in the relation between the 
integer sequence $\bsy{i}$ and the shape $\lambda$ of the tableaux in the image  $DI_n^k(\bsy{i})$. 
\begin{definition}
Let $\bsy{i} \in [n]^k$ be an integer sequence of length $k$, where $n$ and $k$ are positive integers. 
Assume $DI_n^k(\bsy{i})=(P_{\lambda}(\bsy{i}), \Gamma_{\lambda}(\bsy{i}))$ for some integer partition 
$\lambda$. Then we say $\lambda$ is the \emph{VT-shape of $\bsy{i}$}.  
Define the \emph{VT-index of $\bsy{i}$}, denoted by $\mathrm{vt}_n(\bsy{i})$, to be the size of the shape $\lambda^*$. 
That is, $\mathrm{vt}_n(\bsy{i}) $ is the number of boxes in the Young diagram of $\lambda$  that are not in the
first row. \end{definition} 
See Example~\ref{ex:DI} (cf. Section~\ref{Sec:algorithm}) for the integer sequence $\bsy{i}=(3,2,5)$ and its image under $DI_6^3$. 
 For this integer sequence,  the VT-shape is $(3,2,1)$, and hence $\mathrm{vt}_6(\bsy{i}) =3$.

For an integer sequence $\bsy{i}=(i_1,i_2, \ldots, i_k)\in[n]^k$, clearly $0 \leq \mathrm{vt}_n(\bsy{i})\leq k$.  
The exact value of 
the VT-index of $\bsy{i}$ depends on $n$.  
We proved in \cite{BHPYZ1} that when $n \geq 2k + \max(i_1, i_2,\ldots, i_k)$, the vacillating tableau  $\Gamma_{\lambda}^*(\bsy{i})$
is 
stabilized and is called the \emph{limiting \vt}, in which case  $\mathrm{vt}_n(\bsy{i})=k.$

In this paper, we study the cases when $\mathrm{vt}_n(\bsy{i})= 0$ or $k$ for general integers $n$ and $k$. 
Section~\ref{Sec:algorithm} provides the definition of the bijection  $DI_n^k$  as well as its main ingredients:
the RSK row insertion and jeu de taquin.  In Section~\ref{sec: special shapes}, we characterize the sets of integer sequences for three special VT-shapes: 
$\lambda=(n)$, $\lambda=(n-k, 1^k)$ for $n \geq k+1$, and $\lambda=(n-k,k)$ for $n \geq 2k$. 
These sets are in one-to-one correspondences with set partitions of $[k]$ with up to $n$ blocks, the decreasing integer sequences of length $k$ of $[n-1]$, and subdiagonal lattice paths from $(0,0)$ to $(n-k,k)$, respectively. 
Note that the VT-index is $0$ for the first case and is $k$ for the latter two cases. 

In Section~\ref{section: max VT-index}, we characterize  
the set of  
all integer sequences in $[n]^k$ with VT-index $k$ when $n \geq k+1$.      More precisely, for each $\lambda \vdash n$ with $\lambda_1=n-k$,  let $\mathcal{I}_{k}(\lambda)$ be the set of integer sequences in $[n]^k$ whose VT-shape is $\lambda$. 
We introduce a set $\mathcal{R}_k(\lambda)$ of permutations  
and  
show that  an integer sequence 
$\bsy{i} \in [n]^k$  has  VT-shape $\lambda$ 
if and only if  $\bsy{i}$ can be obtained from a permutation in $\mathcal{R}_k(\lambda)$ via some simple transformations.  The transformations are described in Algorithm A of Section~\ref{section: max VT-index}, while its 
inverse is described in Algorithm B. Combining  these two alrgorithms, we develop a test which 
checks whether a sequence $\bsy{i} \in [n]^k$ has  VT-index  $k$. 

An important intermediate object in Algorithms A and B of Section~\ref{section: max VT-index} is the \emph{bumping sequence}. In Section~\ref{sec:bumping}, we characterize the bumping sequences and then reinterpret this characterization in terms of a reparking problem.

Finally, in Section~\ref{sec:final},  we compare our results with another bijective proof of Identity 
\eqref{eq: Identity1}  by Colmenarejo, Orellana,  Saliola, Schilling, and Zabrocki \cite{COSSZ}, which we refer to as the COSSZ bijection. 
 This bijection   maps integer sequences in $[n]^k$ 
to pairs of tableaux of the same shape, where one being an SYT and the other being  a standard multiset tableau.
We explain when an integer sequence corresponds to a shape $\lambda$ with $\lambda_1=n-k$ under the COSSZ  bijection.

\section{The RSK algorithm and the deletion-insertion process} \label{Sec:algorithm}

In this section we recall the RSK algorithm and the bijection $DI_n^k$ constructed by Halverson and Lewandowski \cite{HL05}. 
The  main ingredients of the map $DI_n^k$ are the row insertion algorithm and a special case of  
Sch\"{u}tzenberger's jeu de taquin algorithm, which removes a box containing an entry $x$ from a partial tableau and produces a new partial tableau. 
 We adopt the description from \cite{HL05}, while the full version and in-depth discussion of these algorithms can be found in \cite[Chapter 3]{Sagan01} or \cite[Chapter 7]{EC2}.

\medskip 

\noindent 
\textbf{The RSK row insertion}. Let $T$ be a partial tableau of partition shape $\lambda$ 
with distinct entries. 
Let $x$ be a positive integer that is not in $T$. The operation $x \xrightarrow{RSK}  T$ is defined as follows. 
\begin{enumerate}[(a)]
    \item Let $R$ be the first row of $T$. 
    \item While $x$ is less than some element in $R$:
       \begin{enumerate}[label=\roman*)]
           \item Let $y$ be the smallest element of $R$ greater than $x$;
           
           \item Replace $y \in R$ with $x$;
           \item Let $x:=y$ and let $R$ be the next row. 
       \end{enumerate}
    \item Place $x$ at the end of $R$ (which is possibly empty). 
\end{enumerate}
The result is a partial tableau of shape $\mu$ such that $|\mu /\lambda|=1$. For each occurrence of Step (b), we say that $x$ \emph{bumps} $y$ to the next row.  

\medskip 
The RSK algorithm is an algorithm in algebraic combinatorics that iterates the above insertion procedure. Below is  
Knuth's construction \cite{Knuth70} that applies the algorithm to a two-line 
array of integers  
\begin{equation}  \label{2-line-array}
\left( 
\begin{array}{cccc}
    u_1 & u_2 & \cdots & u_n \\ 
    v_1  & v_2 & \cdots  & v_n 
 \end{array}
\right),
\end{equation}
where $(u_j, v_j)$  are arranged in the non-decreasing lexicographic order from left to right, that is, $u_1 \leq u_2 \leq \cdots \leq u_n$ and $v_j\leq v_{j+1}$ if $u_j=u_{j+1}$. 
\medskip 

\noindent 
\textbf{The RSK algorithm.} \cite{Knuth70} \  
Given the two-line array \eqref{2-line-array}, construct a pair of Young tableaux $(P,Q)$ of the same shape by 
starting with $P=Q=\emptyset$. For $j=1, 2, \ldots, n$:
\begin{enumerate}[(a)]
    \item Insert $v_j$ into tableau $P$ using the RSK row insertion. This 
    operation adds a new box to the shape of $P$. Assume the new box is 
    at the end of the $i$-th row of $P$.
    \item Add a new box with entry $u_j$ at the  end of the $i$-th row  of $Q$.  
\end{enumerate}
We call $P$ the \emph{insertion tableau} and $Q$ the \emph{recording tableau}. 

\medskip
Let $A$ and $B$ be two totally ordered alphabets.   If $u_j \in A$ and $v_j \in B$ for all $j$, we say that the two-line array \eqref{2-line-array} is a generalized permutation from $A$ to $B$.  

\begin{theorem}[Knuth] 
    There is a one-to-one correspondence between generalized permutations from $A$ to $B$ and pairs of SSYTs $(P,Q)$ of the same shape, where content$(P) \subseteq B$ and content$(Q) \subseteq A$.    
\end{theorem}

In particular, there are two special cases of Knuth's RSK algorithm that  correspond to permutations and words, as 
the ones developed by Robinson \cite{Robinson38} and Schensted \cite{Sch61}.  
  
\begin{enumerate}[(i)]
\item When $(u_1\  u_2\ \cdots\ u_n)=(1\  2\ \cdots\ n)$  and $(v_1\  \cdots\  v_n)$ ranges over all permutations of $[n]$,  the correspondence gives a bijection between permutations of length $n$ and pairs of SYTs of the same shape.
\item 
When $(u_1\  u_2 \ \cdots\ u_n)=(1\  2\  \cdots\ n)$ and $v_i \in \mathbb{Z}^+$,  the correspondence gives a bijection between integer sequences of length $n$
and pairs of Young tableaux of the same shape 
$\lambda \vdash n$, where $P$ is an SSYT  with content in $\mathbb{Z}^+$, and $Q$ is an SYT. 
\end{enumerate} 

\medskip

\noindent 
\textbf{Jeu de taquin}. \ 
Let $T$ be a partial tableau of  shape $\lambda$  with distinct entries. Let $x$ be an entry in $T$.  The following operation will delete $x$ from $T$ and yield a partial tableau.
\begin{enumerate}[(a)]
    \item Let $c=T_{i,j}$  be the box of $T$ containing $x$, 
    which is  the $j$-th box from the left in the $i$-th row of $T$. 
    \item While $c$ is not a corner, i.e., a box both at the end of a row and the end of a column, do
      \begin{enumerate}[label=\roman*)]
          \item Let $c'$ be the box containing $\min\{T_{i+1,j}, T_{i, j+1}\}$; if only one of $T_{i+1,j}, T_{i, j+1}$ exists, then the minimum is taken to be that single entry; 
          \item Exchange the positions of $c$ and $c'$. 
      \end{enumerate}
    \item Delete $c$. 
\end{enumerate}
We denote this process by $x \xleftarrow{jdt} T$ and write jdt to denote jeu de taquin.

\medskip

The bijection $DI_n^k$ from the set of integer sequences 
in $[n]^k$ to the set of pairs 
consisting of an SYT and an $n$-vacillating tableau is built by iterating alternatively between the RSK row insertion and jeu de taquin.

\medskip 

\noindent 
\textbf{The bijection $DI_n^k$.} \cite{HL05}  \ 
Let $(i_1, i_2, \ldots, i_k) \in [n]^k$ be an integer sequence of length $k$.  First we define a sequence of partial tableaux recursively: the $0$-th
tableau is the unique SYT of shape $(n)$ with the filling $1, 2, \ldots, n$, namely, 
\begin{center} 
\begin{tikzpicture}
\node[left] at (0,0) {$T^{(0)} =$}; 
\draw (0,-.25) rectangle (3,.25); 
\draw (0.5, -.25)--(0.5,0.25) (1,-0.25)--(1,0.25) (2.5, -0.25)--(2.5, 0.25); 
\node at (0.25,0) {$1$}; 
\node at (0.75, 0) {$2$}; 
\node at (1.75,0) {$\ldots$}; 
\node at (2.75,0) {$n$}; 
\node at (3.2, -0.2) {.}; 
\end{tikzpicture}
\end{center} 

Then for integers $j=0, 1, \ldots, k-1$, the partial tableaux $T^{(j+\frac{1}{2})}$ and $T^{(j+1)}$ are defined by
\begin{eqnarray}
T^{(j+\frac{1}{2})} &= &\left( i_{j+1} \xleftarrow{jdt} T^{(j)} \right), \\ T^{(j+1)}  & = & \left( i_{j+1} \xrightarrow{RSK} T^{(j+\frac{1}{2})} \right). 
\end{eqnarray}
Note that $T^{(j+1)}$ is always an SYT. 
Let $\lambda^{(m)}$ be the shape of $T^{(m)}$ for all integral and half-integral indices $m$,   and $\lambda=\lambda^{(k)}$ be the shape of the last partial tableau $T^{(k)}$. 
Finally, let 
\begin{eqnarray}
P_{\lambda}(\bsy{i})=T^{(k)}, \qquad \Gamma_{\lambda}(\bsy{i})
=\left(  
\lambda^{(0)}, \lambda^{(\frac{1}{2})}, \lambda^{(1)}, \lambda^{(1\frac{1}{2})}, \ldots, \lambda^{(k)} 
\right),  
\end{eqnarray}
so $P_{\lambda}(\bsy{i})$ is an SYT of shape $\lambda$ and 
$\Gamma_{\lambda}(\bsy{i})$ is an $n$-vacillating tableau of shape $\lambda$ and length $2k$.  The image of $\bsy{i}$ under the map $DI_n^k$ is given by  $DI_n^k(\bsy{i}) =  \left(P_{\lambda}(\bsy{i}), 
\Gamma_{\lambda}(\bsy{i})\right)$. 

Note that the map $DI_n^k$  is a well-defined bijection for all positive integers $n$ and $k$. Hence, Identity~\eqref{eq: Identity1} is indeed true for all $n, k \in \mathbb{Z}_{>0}$.

\begin{example}  \label{ex:DI} 
Let  $n = 6$, $k=3$  and $\bsy{i}=(3,2,5)$.  We apply the map $DI^3_6$ to obtain the $(T^{(j)})$ sequence as follows.
 \[
 \ytableausetup{smalltableaux}
 \begin{ytableau}
1 &2 &3 &4 &5 & 6 
\end{ytableau} \ \to\ \begin{ytableau}
1 &2 &4 &5 & 6
\end{ytableau}\ \to\  \begin{ytableau}
1 &2 &3 &5 & 6 \\
4
\end{ytableau}\ \to\  \begin{ytableau}
1  &3  & 5 & 6\\
4
\end{ytableau}\ \to\ \begin{ytableau}
1 &2 &5 &6\\
3 \\
4 
\end{ytableau} \ \to\ \begin{ytableau}
1 &2  &6\\
3 \\
4 
\end{ytableau} \  \to\ \begin{ytableau}
1 &2  &5\\
3 & 6 \\
4 
\end{ytableau}  
\] 
Thus, for $DI_6^3(\bsy{i})$,  $\lambda=(3,2,1)$,  and 
\[ 
P_{\lambda}(\boldsymbol{i})= \ytableausetup{smalltableaux} \begin{ytableau}
1 &2 &5  \\
3 & 6\\
4
\end{ytableau},  \quad  
 \Gamma_{\lambda}(\boldsymbol{i})=  \ytableausetup{boxsize=7pt}  \left(\ydiagram{6},\, \ydiagram{5}, \ \ydiagram{5,1},\ \ydiagram{4,1},\ \ydiagram{4,1,1},\ \ydiagram{3,1,1}, \ \ydiagram{3,2,1} \right). \] 
\end{example}

\section{Integer sequences with special VT-shapes}\label{sec: special shapes}

Let $n$ and $k$ be positive integers.  Given a partition $\lambda \vdash n$, we are interested in characterizing integer sequences $\bsy{i} \in [n]^k$ with  the VT-shape $\lambda$.  Let $\mathcal{I}_{k}(\lambda)$ be the set of all such integer sequences~$\bsy{i}$. 
From the definition of $DI_n^k$, we have  that $DI_n^k$ induces a bijection between $\mathcal{I}_{k}(\lambda)$ 
and $\mathcal{SYT}(\lambda) \times \mathcal{VT}_{n,k}(\lambda)$, where $\mathcal{SYT}(\lambda)$ is the set of all SYT with shape $\lambda$. As a result,
\begin{equation}\label{eq:size of I_k}
|\mathcal{I}_{k}(\lambda)|= f^\lambda \cdot |\mathcal{VT}_{n,k}(\lambda)|.
\end{equation}

Next we characterize the set $\mathcal{I}_k(\lambda)$ for three special shapes: 
$\lambda =(n), (n-k,1^k)$, and $(n-k,k)$.

\noindent 
\emph{Case 1. $\lambda =(n)$. } 

It is proved in \cite{BHH17,MR98} that the number $|\mathcal{VT}_{n,k}((n))|$ of $n$-vacillating tableaux with length $2k$ and shape $(n)$ equals $\sum_{i \leq n} S(k,i)$, where $S(k,i)$ is the Stirling number of the second kind and counts the number of set partitions of $[k]$ into exactly $i$ blocks.  
When $n \geq k$, the sum $\sum_{i \leq n} S(k,i)$ is the Bell number $B(k)$, and the equation 
$|\mathcal{VT}_{n,k}((n))|=B(k)$ is proved in  \cite[Theorem~2.4]{CDDSY07}.
Since $f^{(n)} = 1$, by Equation \eqref{eq:size of I_k}, we have 
\[
|\mathcal{I}_k((n))|= f^{(n)}\cdot |\mathcal{VT}_{n,k}((n))| =  \sum_{i \leq n} S(k,i).
\] 
Theorem~\ref{theorem:one-row}  characterizes all the  sequences $\bsy{i} \in \mathcal{I}_k((n))$.

\begin{theorem}\label{theorem:one-row}
Let $n$ and $k$ be positive integers.  An integer sequence $\bsy{i}=(i_1, i_2, \ldots, i_k) \in [n]^k$ is in $\mathcal{I}_k((n))$ if and only if $\bsy{i}$ satisfies the following condition: For each $r$ such that $1\leq r \leq k$, if $i_r = m$ for some $m < n$,  then there exists an integer $s$ such that $r < s \leq k$ and $i_s = m+1$. In particular, $i_k = n$. 
\end{theorem}
\begin{proof} \ 
 We first show that if $\bsy{i}\in \mathcal{I}_k((n))$, then it must have the described property. 
 It is easy to check that 
 the sequence $(n,\ldots, n)$ is in $\mathcal{I}_k((n))$ and satisfies the property. If $\bsy{i} \neq (n,\ldots, n)$,  let $m <n$ 
be an integer appearing in $\bsy{i}$ and $r \leq k $ be the latest position such that $i_r = m$.  
Then in the delete-insert process of $DI_n^k$ when $i_r=m$ is removed by jdt and then inserted back using the RSK row insertion, 
we obtain the tableau $T^{(r)}$, in which $m$ is in the first row and $m+1$ is in a row strictly below $m$. 
Assume that the integer $m+1$ does not appear in $i_{r+1}, \ldots, i_k$. 
Then we can prove by induction   that in the remaining rounds of deletion-insertion
 $m+1$ always stays below $m$.  For $j=r+1, \ldots, k$:
\begin{itemize}
   \item When deleting $i_j$ by jdt from $T^{(j)}$ : 
       \begin{enumerate}[(i)]
        \item Assume $i_j< m$. The operation of jdt does not move any entry to a row with a larger index. If in the process of jdt, $m+1$ is moved up by exchanging with $i_j$ from position $(t+1, s)$ to $(t, s)$, then we claim that $m$ cannot be in row $t$. Otherwise, assume $m$ is at $(t, x)$ with $x < s$. Then the box $(t+1, x)$ is in the 
        tableau $T^{(j)}$ and occupied by an entry $y$ such that $ m< y < m+1$, which is impossible. Therefore $m$ must be in a row with index less than $t$. 

        If in the process of jdt, $m+1$ is not moved up, then it stays below $m$.

        \item  Assume $i_j>m+1$. Then $i_j$ is only exchanged with entries $c > i_j$. Hence the process does not change the position of $m$ or $m+1$. 
        \end{enumerate} 
    \item  When inserting $i_j$ to $T^{(j+\frac12)}$ by row insertion:   
  \begin{enumerate}[(i)]
    
      \item  Assume $i_j < m$.  If the insertion path of $i_j$ contains $m$, assume $m$ is bumped from row $t$ to row $t+1$. Then  $m+1$ is either already at a row below $t+1$, or bumped by $m$ to row $t+2$.  

       If the insertion of  $i_j$ does not contain $m$, since the row insertion algorithm will not move any entry to a row with smaller index, $m+1$ remains below $m$. 

      \item Assume $i_j>m+1$.  Then $i_j$ only bumps entries larger than $m+1$ and does not change the positions of  $m$ and  $m+1$.  
    \end{enumerate} 
    \end{itemize}
In any case, $m+1$ stays below $m$.
Therefore we cannot get the shape $(n)$ at the end, which is a contradiction.  

Conversely, if $\bsy{i}$ is a sequence satisfying the property in the statement, for any integer $m$ appearing in $\bsy{i}$ whose last appearance is at position $r$, delete $m$ and insert $m$ will leave $m$ in the first row of $T^{(r)}$.  Since all the later integers are larger than $m$, $m$ will stay at the same box of the first row until the end of the process. The same is true for any integer $t$ that is not appearing in $\bsy{i}$: it will stay in the $t$-th box of the first row throughout the process. Therefore the ending tableau $T^{(k)}$ must have shape $(n).$ 
\end{proof}

Notice that a sequence satisfying the condition in Theorem~\ref{theorem:one-row}
corresponds uniquely to a set partition of $[k]$, where $i_r = i_s$ if and only if $r$ and $s$ are in the same block. 
If we order the blocks in a set partition according to their largest elements in decreasing order, then the correspondence can be written as: $i_r = m$ if and only if $r$ belongs to the $(n-m+1)$-th block. 
Hence, $\mathcal{I}_k((n))$ are in one-to-one correspondence with set partitions of $[k]$ with up to $n$ blocks, which gives $| \mathcal{I}_n^k((n))| = \sum_{i \leq n} S(k,i)$.

\medskip 

\noindent 
\emph{Case 2. $\lambda = (n-k, 1^k)$ with $n \geq k+1$, i.e., $\lambda$ is a hook shape with $k+1$ cells in the first column. } 

There is exactly one $n$-vacillating tableau of length $2k$ with shape $(n-k, 1^k)$, namely, $\lambda^{(i)}=(n-i, 1^i)$ and $\lambda^{(i+\frac12)}=(n-i-1, 1^i)$ 
for $i=0, 1, \ldots, k-1$. 
Hence, $|\mathcal{VT}_{n,k}((n-k, 1^k))| =1$, which by Equation \eqref{eq:size of I_k} implies  \[|\mathcal{I}_k((n-k, 1^k))| = f^{(n-k, 1^k)} = \binom{n-1}{k}.
\]

\begin{theorem}\label{theorem:hook}
For positive integers $n$ and $k$ such that $n \geq k+1$, an integer sequence $\bsy{i}=(i_1, i_2, \ldots, i_k) \in \mathcal{I}_k((n-k,1^k))$ if and only if $n>i_1>i_2>\cdots > i_k.$ 
\end{theorem}
\begin{proof} 
Let $\bsy{i} = (i_1,\ldots,i_k)$ be such that $n > i_1 > \cdots > i_k$. Let $\Gamma(\bsy{i}) = (T^{(0)},T^{(\frac{1}{2})},T^{(1)},T^{(1\frac{1}{2})},\ldots,T^{(k)})$ denote the sequence of partial tableaux obtained by applying $DI_n^k$ to $\bsy{i}$. Since $i_1 < n$, removing $i_1$ by jdt and inserting it back by RSK will make $T^{(1)}$ an SYT of shape $(n-1,1)$ with $i_1+1$ at position $(2,1).$ We assume, for induction, that $T^{(j)}$ is an SYT of shape $(n-j,1^{j})$ with $i_{j}+1$ at position $(2,1)$. Since $i_{j+1} < i_{j+1} +1 < i_j+1$, by removing $i_{j+1}$ from $T^{(j)}$ using jdt, we  obtain a partial tableau of shape $(n-j-1,1^{j})$ that still has $i_{j} + 1$ at position $(2,1)$. Inserting $i_{j+1}$ into $T^{(j+\frac{1}{2})}$ using RSK bumps $i_{j+1} + 1$ from the first row to position $(2,1)$ and bumps each integer at position $(s,1)$ to $(s+1,1)$ for all $2\le s \le j$. 
Hence, $T^{(j+1)}$ is of shape $(n-j-1,1^{j+1})$ with $i_{j+1}+1$ at position $(2,1)$ and this completes the induction. 
It follows that $T^{(k)}$ has the shape $(n-k, 1^k)$.

Let 
$I = \{(i_1, \ldots, i_k) \mid n > i_1 > \cdots > i_k\}$. 
The above paragraph shows that $I \subseteq \mathcal{I}_k((n-k,1^k))$. But $|I|=\binom{n-1}{k} = 
| \mathcal{I}_k((n-k,1^k))|$. It follows that $I = \mathcal{I}_k((n-k,1^k))$, as desired.
\end{proof}

\noindent 
\emph{Case 3. $\lambda = (n-k, k)$ with $n \geq 2k$.  } 

There is exactly one $n$-vacillating tableau  of length $2k$ with shape $(n-k, k)$, namely, 
$\lambda^{(i)}=(n-i, i)$ and $\lambda^{(i+\frac12)}=(n-1-i, i)$ for $i=0, 1, \ldots, k-1$. Therefore, $|\mathcal{VT}_{n,k}((n-k, k))| =1$ and $|\mathcal{I}_k((n-k, k))| = f^{(n-k, k)}$. 

Given an SYT $P$ of shape $(n-k,k)$ with the second row  being $\bsy{b}=(b_1, b_2, \ldots, b_k)$, 
we obtain a sequence $\bsy{i}$  by dividing $(b_1, b_2, \ldots, b_k)$ into disjoint maximal contiguous segments and replace each segment $a, a+1, \ldots, a+\ell-1$ with $a-1, a-1, \ldots, a-1$ ($\ell $ copies). We show in Theorem~\ref{thm:(n-k,k)} that the sequence $\bsy{i}$ obtained by the above operation is in $\mathcal{I}_k((n-k,k))$.  
In fact, $P$ is the SYT component $P_\lambda(\bsy{i})$ when we apply $DI_{n}^{k}$ to $\bsy{i}$.

\begin{example}\label{SYT1} Let $n = 15$ and $k = 7$. Consider the SYT \[P = 
\ytableausetup{aligntableaux=center,boxsize=.5cm}\begin{ytableau}
1 & 2 & 4 & 7 & 8 & 9 & 11 & 15\\
3 & 5 & 6 & 10 & 12 & 13 & 14
\end{ytableau}\]
whose  second row is $\bsy{b}=(3, 5, 6, 10, 12, 13, 14)$. The segments are
$3 \,| \, 5, 6\, |\, 10 \,|\, 12, 13, 14$ which gives $\bsy{i}= (2, 4, 4, 9, 11, 11, 11)$. 
One can check that in $DI_{15}^{7}(\bsy{i})$, the SYT is exactly $P$.
\end{example}

\begin{theorem}  \label{thm:(n-k,k)}
The above operation is a  bijection between $\mathcal{I}_k((n-k, k))$ and the set of SYTs of shape $(n-k,k)$. 
\end{theorem}

\begin{proof}[Proof of Theorem~\ref{thm:(n-k,k)}]
Notice that any SYT of shape $(n-k,k)$ can be represented by a lattice path from $(0,0)$ to $(n-k,k)$ that stays on or below the line $y = x$. The lattice path is recorded by a sequence of east ($E$) and north ($N$) steps, say, $P_1P_2\cdots P_n$, where $P_i = E$ if $i$ is in the first row of the SYT and $P_i = N$ if $i$ is in the second row. Let $\bsy{v} = (v_1,v_2,\ldots, v_k)$ be the $x$-coordinates of the $N$-steps in the lattice path. Then $\bsy{v}$ is a weakly increasing sequence larger than or equal to $(1,2,\ldots,k)$ coordinate-wise. The vectors $\bsy{b}$ and $\bsy{v}$ determine each other by the simple relation that $b_i = i+v_i.$  We write \begin{equation} \label{vector-v}
\bsy{v}=a_1^{r_1} a_2^{r_2} \cdots a_s^{r_s},  
\end{equation} where $a_1<a_2<\cdots <a_s$. Note that we must have 
$a_1 \geq r_1, a_2 \geq r_1+r_2, \ldots, a_s \geq r_1 + \cdots +r_s$.  Then the sequence $\bsy{i} $ obtained by dividing $(b_1,b_2,\ldots,b_k)$ into disjoint maximal contiguous segments and replacing each segment $a, a+1, \ldots, a+\ell-1$ with $a-1, a-1, \ldots, a-1$ ($\ell $ copies) can be written as  
\[
\bsy{i} = \bsy{v}+ \bsy{\epsilon}, 
\]
where 
\[
\bsy{\epsilon} = 0^{r_1} r_1^{r_2} (r_1+r_2)^{r_3} \cdots (r_1+ \cdots + r_{s-1})^{r_s}. 
\]

\noindent 
\textbf{Claim}: For the above $\bsy{i}$, $\bsy{v} $ and $\bsy{b}$, applying $DI_n^k$ to $\bsy{i}$, we have 
$P_\lambda(\bsy{i}) = P$, where $P$ is the SYT of shape $\lambda=(n-k,k)$ whose second row is $\bsy{b}$.

\noindent \textit{Proof of the Claim.} If $\bsy{v}$ is given by \eqref{vector-v}, then 
\[
\bsy{i} = a_1^{r_1} (a_2+r_1)^{r_2} \cdots (a_s+ r_1+ \cdots+r_{s-1})^{r_s}
\]
and \[
\bsy{b} = (a_1+1, \ldots, a_1+r_1, a_2+r_1+1, \ldots, a_2+r_1+r_2, \ldots, 
a_s+1+r_1+\cdots+r_{s-1}, \ldots, a_s+r_1+\cdots +r_s).
\]

Let $\Gamma(\bsy{i}) = (T^{(0)},T^{(\frac{1}{2})},T^{(1)},T^{(1\frac{1}{2})},\ldots,T^{(k)})$ be the sequence of partial tableaux obtained by applying $DI_n^k$ to $\bsy{i}$. 
We  show by induction that for each integer $i$, $T^{(i)}$ is the SYT of shape $(n-i, i)$ whose second row consists of the first $i$ entries of $\bsy{b}$. Since $a_1 < n$, removing $a_1$ by jdt and inserting it back by RSK results in $T^{(1)}$ being an SYT of shape $(n-1,1)$ with $a_1+1$ at position $(2,1)$. Let $i-1 \ge 1$. Assume that $T^{(i-1)}$ is the SYT of shape $(n-i+1,i-1)$ whose second row comprises the first $(i-1)$ entries of $\bsy{b}$. We can write $i = r_1+\cdots +r_{j-1} +m$ for some $j\ge 1$ and $1\le m\le r_j$. The $i$-th entry of $\bsy{i}$ is $a_j+r_1+\cdots + r_{j-1}$.

Now consider $T^{(i)}$. By the inductive hypothesis, the largest entry in the second row of $T^{(i-1)}$ is $a_j+i-1 = a_j+r_1+\cdots +r_{j-1} +m-1$. Since there are $i -1$ entries in the second row and the largest entry in the second row is $a_j + i -1$,  there are at least $(a_j+i-1)-(i-1) = a_j \ge r_1+r_2+\cdots +r_j \ge i$ entries less than $a_j+i-1$ in the first row of $T^{(i-1)}$.

As $a_j+r_1+\cdots + r_{j-1}$ is the largest integer less than $a_j+i-1$ in the first row, there are at least $i-1$ entries less than $a_j+r_1+\cdots + r_{j-1}$ in the first row of $T^{(i-1)}$. Consequently, $a_j+r_1+\cdots +r_{j-1}$ is in the first row with no entry right below it. Therefore, removing $a_j+r_1+\cdots+r_{j-1}$ by jdt and inserting it back by RSK will move the smallest entry greater than $a_j+r_1+\cdots + r_{j-1}$ in the first row, i.e., $(a_j + i -1) + 1  =a_j+r_1+\cdots+r_{j-1}+m$, to the second row, resulting in the desired pattern in $\bsy{b}$.
\end{proof}

We illustrate the proof of Theorem~\ref{thm:(n-k,k)} next.

\begin{example} \label{ex:lattice} 
For the integer sequence in Example~\ref{SYT1}, $\bsy{i}= (2, 4, 4, 9, 11, 11, 11) = 2^14^29^111^3$ and hence $\bsy{\epsilon} = 0^1 1^2 3^1 4^3$.  The $x$-coordinates of the North steps are $\bsy{v} = (2,3,3,6,7,7,7)$. We illustrate the corresponding lattice path in Figure~\ref{fig:example lattice}. We can write $\bsy{v} =2^13^26^17^3$. Then the sequence $\bsy{i}= \bsy{v} +\bsy{\epsilon} = (2, 4, 4, 9, 11, 11, 11)$, as expected.
One can check that the lattice path determined by $\bsy{v}$ corresponding to SYT $P = P_{(8,7)}(\bsy{i})$ in Example~\ref{SYT1}. 

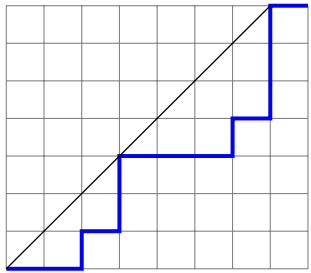
\begin{figure}[h]
    \centering
    \begin{tikzpicture}[scale = 0.5,anchor=base, baseline]
 \draw[step=1cm, color=gray] (0, 0) grid (8,7);
\draw[line width=1.5pt, color=blue] (0,0) -- (2,0) -- (2,1) -- (3,1) -- (3,3) -- (6,3) -- (6,4) -- (7,4) -- (7,7) -- (8,7);
\draw[line width=.5pt, color=black] (0,0)-- (7,7);
\end{tikzpicture}
    \caption{Lattice path corresponding to $\bsy{v} = (2,3,3,6,7,7,7)$.}
    \label{fig:example lattice}
\end{figure}

\end{example}

The proof of Theorem~\ref{thm:(n-k,k)} also establishes a correspondence between sequences in $\mathcal{I}_k((n-k,k))$ and lattice paths from $(0,0)$ to $(n-k,k)$ below the line $y=x$. Let $\bsy{i}=(i_1, i_2, \ldots, i_k) \in \mathcal{I}_k((n-k,k))$. Assume  $\bsy{i}=c_1^{r_1} c_2^{r_2} \cdots c_s^{r_s}$, where $c_1 < c_2 < \cdots < c_s$. Let  $\bsy{\epsilon} = 0^{r_1} r_1^{r_2} (r_1+r_2)^{r_3} \cdots (r_1+ \cdots + r_{s-1})^{r_s}$. Then, $i_j - \epsilon_j$ is the $x$-coordinate of the North step of a lattice path from $y= j-1$ to $y=j$ for $1\leq j \leq k$. 
It follows that $|\mathcal{I}_k((n-k,k))| = \frac{n-2k-1}{n-k+1} \binom{n}{k}$, \cite{Handbook10}. 
In particular,  when  $n = 2m$ and $k = m$,  this gives a bijection between $\mathcal{I}_k((m,m))$  and  the set of Dyck paths from $(0,0)$ to $(m,m)$, which is counted by the $m$-th Catalan number.

\section{Integer sequences with maximal VT-index}\label{section: max VT-index}
In Section~\ref{sec: special shapes},  we characterize all $\bsy{i} \in [n]^k$ with VT-index $0$ and two special cases when 
the VT-shape is a hook or a two-row partition,  both of which have VT-index $k$. In this section, we characterize  all integer sequences in $[n]^k$ with VT-index $k$.   More precisely, 
for each $\lambda \vdash n$ with $\lambda_1=n-k$,  (or equivalently, $|\lambda^*|=k$),  
we introduce a set $\mathcal{R}_k(\lambda)$ of permutations, which has the same cardinality as $\mathcal{I}_k(\lambda)$. 
We establish a bijection between these two sets and show that for $n \geq k + 1$, an integer sequence 
$\bsy{i} \in [n]^k$  has  the VT-shape $\lambda$ 
if and only if  $\bsy{i}$ can be obtained from a permutation in $\mathcal{R}_k(\lambda)$ via some simple transformations, as defined 
in Algorithm A in this section.

\begin{lemma}\label{lemma:size_vt}
Let $\lambda^*=(\lambda_2, \lambda_3,\ldots, \lambda_l)$ be an integer partition of a positive integer $k$.  Let $n \geq k+\lambda_2$ and $\lambda=(n-k, \lambda_2,\lambda_3, \ldots, \lambda_l)$ be an integer partition of  $ n$.  
Then the number of $n$-vacillating tableaux of shape $\lambda$ and length $2k$ is $f^{\lambda^*}$, 
the number of SYTs of shape $\lambda^*$. Consequently, 
 $|\mathcal{I}_k(\lambda)| = f^\lambda \cdot f^{\lambda^*}$.
\end{lemma}
\begin{proof}
    Given  $\Gamma_\lambda =  ( \lambda^{(j)}: j=0, \frac{1}{2}, 1, 1\frac{1}{2}, \dots, k)$ in $\mathcal{VT}_{n,k}(\lambda)$,  let $\Gamma_\lambda^*$ be the corresponding simplified
    \vt  \
     $( \mu^{(j)}: \mu^{(j)}= (\lambda^{(j)})^* \text{ for  } j=0, \frac{1}{2}, 1, 1\frac{1}{2}, \dots, k)$. 
    Since $\vert \mu^{(k)}\vert = \vert \lambda^* \vert =k$, we must have that 
    $\mu^{(j)}= \mu^{(j+\frac12)}$ and $\vert \mu^{(j+1)}/\mu^{(j)}\vert =1$ for all $j=0, 1, \dots, k-1$. 
    Hence $\Gamma_\lambda$ is in one-to-one correspondence with  the SYT $Q^*$ of $\mu^{(k)}=\lambda^*$, where the integer $j \in [k]$ is at the box of $\mu^{(j+1)}/\mu^{(j)}$ in $Q^*$.    This gives $|\mathcal{VT}_{n,k}(\lambda)| = f^{\lambda^*}$.
\end{proof}

Let $n \geq k+1$. For a permutation $w =w_1w_2\cdots w_n\in \mathfrak{S}_n$ written in the one-line notation, 
let $\mathrm{is}(w)$ be the length of the longest increasing subsequences of $w$. Define $\mathcal{R}_{k}^n \subseteq \mathfrak{S}_n$ by letting   
\[
\mathcal{R}^n_{k}=\{ w_1 w_2\cdots w_n  \in \mathfrak{S}_n: \  w_1 < w_2 < \cdots < w_{n-k} \text{ and } \mathrm{is}(w)=n-k\}.  
\]
Let $(P(w), Q(w))$ be the pair of SYTs obtained by applying the RSK algorithm to $w$, where $P(w)$ is the insertion tableau and $Q(w)$ is the recording tableau. Then $w \in \mathcal{R}^n_{k}$ if and only if $\lambda$, the shape of $P(w)$ and $Q(w)$, 
satisfies $\lambda_1=n-k$ and the entries in the first row of $Q(w)$ are $1, 2, \ldots, n-k$. Let $\mathcal{R}_k(\lambda)$
be the permutations $w \in \mathcal{R}^n_{k}$ such that the shape of $P(w)$ is $\lambda$. Then we have 
$|\mathcal{R}_k(\lambda)| = f^\lambda \cdot f^{\lambda^*}$ and 
\[
\mathcal{R}^n_{k}= \biguplus_{\lambda \vdash n,\, \lambda_1=n-k} \mathcal{R}_k(\lambda). 
\]

Note that for $\lambda \vdash n$ with $\lambda_1 =n- k$, $\mathcal{I}_k(\lambda)$ and $\mathcal{R}_k(\lambda)$ have the same size.  We  define an explicit  bijection $\psi$ between these two sets, and then use the permutations in $\mathcal{R}_k(\lambda)$ to 
characterize the integer squences in $\mathcal{I}_k(\lambda)$.   

In the rest of this section, we always assume that $\lambda \vdash n$ with $\lambda_1 =n- k$, (and hence $|\lambda^*|=k$). 
\vspace{.2cm} 

\noindent \textbf{A bijection $\psi$ from $\mathcal{I}_k(\lambda)$ to $\mathcal{R}_k(\lambda)$.}

Let $\bsy{i} \in \mathcal{I}_k(\lambda)$.
\begin{enumerate}[(a)]
    \item  Let
$DI_n^k(\bsy{i})=(P_\lambda(\bsy{i}), \Gamma_\lambda(\bsy{i}))$,  where $P_\lambda(\bsy{i}) \in \mathcal{SYT}(\lambda)$ and $\Gamma_\lambda(\bsy{i}) \in \mathcal{VT}_{n,k}(\lambda).$

\noindent 
Assume $\Gamma=\Gamma_\lambda(\bsy{i})=(\lambda^{(i)}: i=0, \frac{1}{2}, 1, \ldots, k)$.  Let 
$\Gamma^*=\Gamma^*_\lambda(\bsy{i})= (\mu^{(i)}=(\lambda^{(i)})^*: i=0, \frac{1}{2}, 1, \ldots, k)$ be the corresponding simplified
\vt. 
Then $\mu^{(0)}=\emptyset$, $\mu^{(k)}=\lambda^*$, and for each integer $i=0, 1, \ldots, k-1$, 
$\mu^{(i)} \subset \mu^{(i+1)}$ and $|\mu^{(i)}|=i$. 
\item 
Let $Q^*$ be an SYT of shape $\lambda^*$ that corresponds to the sequence $\mu^{(0)}, 
\mu^{(1)}, \ldots, \mu^{(k)}$. That is, entry $i$ of $Q^*$ is in the unique box of
$\mu^{(i)} / \mu^{(i-1)}$. 
\item Create a SYT $Q$ of shape $\lambda$ from $Q^*$ by
adding $n-k$ to all the $k$ entries  in $Q^*$ (so all these entries form the set $\{n-k+1,n-k+2, \ldots, n\}$), followed by adding a top row with entries $1,2, \ldots, n - k$.
\item  Let $w$ be the unique permutation such that $RSK(w)=(P_\lambda(\bsy{i}), Q)$. Then $w \in \mathcal{R}_k(\lambda)$. 
\item We define $\psi(\bsy{i})=w$. 
\end{enumerate}

\begin{example} \label{ex:psi}
    Let $n=8$, $k=4$, and  $\bsy{i}=(3,2,6,2)$. Then $DI_8^4(\bsy{i})=(P_\lambda(\bsy{i}), \Gamma_\lambda(\bsy{i}))$  of shape $\lambda=(4,2,2)$ such that 
    \[
 P_\lambda(\bsy{i})=  \ytableaushort{1268,35,47}
        \]
        and 
        \[
\Gamma_\lambda(\bsy{i})=
\left(\ytableausetup{boxsize=6px} 
 \ydiagram{8}\, ,\, \ydiagram{7}\, ,\, \ydiagram{7,1}\, ,\, \ydiagram{6,1}\, ,\, \ydiagram{6,1,1}\, ,\, \ydiagram{5,1,1}\, ,\, \ydiagram{5,2,1}\, ,\, \ydiagram{4,2,1}\, ,\, \ydiagram{4,2,2}\ \right).
\]
Therefore 
\[ \ytableausetup{boxsize=12px} 
Q^*= \ytableaushort{13,24} \qquad \text{ and } \qquad 
Q=  \ytableaushort{1234,57,68}.
        \]
 The pair $(P_\lambda(\bsy{i}), Q)$ corresponds to $w=45783162$ under RSK.     Hence $\psi( (3,2,6,2))= 45783162$.     
\end{example}

\noindent \textbf{The inverse of $\psi$. }\ 
The map $\psi$ is a bijection since all the steps can be reversed. Starting with $w \in \mathcal{R}_k(\lambda)$, 
we can recover the integer sequence $\bsy{i}$ via the following steps. 
\begin{enumerate}[(a$'$)]
    \item Apply the RSK algorithm to $w$ to get the insertion tableau $P(w)$ and the recording tableau $Q(w)$. 
    Let  $\lambda$ be the shape of $P(w)$. 
    \item Note that the first row of $Q(w)$ must be exactly $1, 2, \ldots, n-k$. Removing the first row and reducing every remaining entry of $Q(w)$ by $n-k$, we get an SYT $Q^*$ of shape $\lambda^*$. 
    \item Let $\Gamma'= (\mu^{(i)}: i=0, \frac{1}{2}, 1, \ldots, k)$ be the simplified \vt\ of length $2k$ and shape
    $\lambda^*$ by letting $\mu^{(k)}=\lambda^*$, $\mu^{(0)}=\emptyset$,  and $\mu^{(i)}= \mu^{(i +\frac12)}$ be the shape of $Q^*$
    restricted to the boxes with entries $\{1,2,  \ldots, i\}$, for $i=1, 2, \ldots, k-1$. 
    \item Let $\Gamma$ be the $n$-\vt\ whose simplified version is $\Gamma'$.  
    \item Then $\bsy{i}=\psi^{-1}(w)$ is  be the  sequence such that $DI_n^k(\bsy{i})=(P(w), \Gamma)$. 
    Since $DI_n^k$ is a bijection, $\bsy{i} \in [n]^k$ exists and is unique.  

\end{enumerate}

The bijection $\psi$ gives a one-to-one correspondence between $\mathcal{I}_k(\lambda)$ and 
$\mathcal{R}_k(\lambda) \subseteq\mathfrak{S}_n$. The above definition uses the map $DI_n^k$
and the RSK algorithm, hence it is not easy to obtain a characterization of  the integer sequences
in $\mathcal{I}_k(\lambda)$. 

Next we give a more direct description of the correspondence $\psi$. We show that
the sequence $\bsy{i} \in \mathcal{I}_k(\lambda)$ and $w=\psi(\bsy{i}) \in \mathcal{R}_k(\lambda)$ can be obtained  from each other
by some simple operations without using  jeu de taquin or  row insertion.  
 Algorithm A below provides the operations transforming  $w$ to $\bsy{i}$, while Algorithm B is the reversed
operations of Algorithm A. A key ingredient of these algorithms is an intermediate sequence 
$\bsy{t}=(t_1, t_2, \ldots, t_k)$, which is called \emph{a bumping sequence} because they are the integers that get bumped out from the first row in the last $k$ rounds when we apply the RSK algorithm in the process of $w \xrightarrow[]{RSK} (P(w), Q(w))$.  
From the correspondence $\psi$, $t_j$ is also the integer that gets bumped out from the first row when we apply 
the row-insertion $i_j \rightarrow T^{(j-\frac12)}$ in the process of $DI_n^k(\bsy{i})$ for $\bsy{i} \in \mathcal{I}_k(\lambda)$.

\begin{center}
\fbox{\begin{minipage}[t][1.1\height][c]{\dimexpr\textwidth-2\fboxsep-2\fboxrule\relax}
\noindent \textbf{Algorithm A: from $\mathcal{R}_k(\lambda)$ to $\mathcal{I}_k(\lambda)$. } 
\begin{enumerate}
    \item Start with   $w = b_1 b_2\cdots b_{n-k} a_1  a_2\cdots a_k  \in \mathcal{R}_k(\lambda)$.
    
    \item For $j=1$ to $k$, let 
    \[
     t_j := \min\{ x \in \mathbb{Z}: \ x > a_j \text{ and } x \neq t_1, \ldots, t_{j-1}, a_{j+1}, \ldots, a_k\} .
    \] 
    \item For $j=1$ to $k$, let 
    \[
    i_j := \max\{ c \in \mathbb{Z}:\  c < t_j \text{ and } c \neq t_1, t_2, \ldots, t_{j-1} \}. 
    \]
    \item The output is the integer sequence $\bsy{i}=(i_1, i_2,\ldots, i_k)$.
\end{enumerate}
\end{minipage}}
 \end{center}

The above procedure can be  visualized using an $n \times k$ rectangular grid. Let $(r,c)$ be the cell in 
row $r$ from the bottom and column $c$ from the left. 
Given  $w = b_1\cdots b_{n-k}\ a_1\cdots a_k \in \mathcal{R}_k(\lambda)$,  
we use diagonal shading 
\raisebox{-3pt}{\tikz{\fill[pattern=north east lines, pattern color=lightgray] (0,-1) rectangle (0.5,-1.5);\draw(0,-1) rectangle (0.5,-1.5)}} to shade the cell $(a_j, j)$  and all the cells in the same row on 
its left $\{(a_j,j') \ | \  1 \leq j' < j \leq k\}$ for all $1 \leq j \leq k$. 
Starting from $j = 1$ to $j = k$, let $t_j$ be the row index (from bottom) of the lowest un-shaded cell above 
$(a_j,j)$ in column $j$, and 
use horizontal shading 
\raisebox{-3pt}{\tikz{\fill[pattern=horizontal lines light gray, pattern color=lightgray] (0,-1) rectangle (0.5,-1.5);\draw(0,-1) rectangle (0.5,-1.5)}} to
shade the cells on the right of $(t_j,j)$, that is, cells $(t_j, j')$ for all 
$j \leq j' \leq k)$.  
After all the cells $(t_j,j)$ are determined, remove the diagonal shading but  keep the horizontal ones, and 
 color in black \raisebox{-3pt}{\tikz{\fill[black] (0,-1) rectangle (0.5,-1.5);\draw(0,-1) rectangle (0.5,-1.5)}} the highest non-shaded cell $(i_j,j)$ below cell $(t_j,j)$. 
 Then $\bsy{i}=(i_1, i_2, \ldots, i_k) \in \mathcal{I}_k(\lambda)$.

\begin{example}
The following figures demonstrate the implementation of Algorithm A.  Let $n=8$, $k=4$ and 
$w = 45783162\in \mathcal{R}_4((4,2,2))$. We start with $(a_1, a_2, a_3, a_4)=(3, 1, 6, 2)$. 
The process yields $\bsy{t} = (4,3,7,5)$ and then 
 $\bsy{i} = (3,2,6,2) \in \mathcal{I}_4((4,2,2))$, which corresponds to the black cells in the last grid.

\input{image_revised}

\input{keyforfigures}

\end{example}

\begin{theorem}
   Let $\lambda \vdash n$ with $\lambda_1 =n-k$. 
    Assume that $\bsy{i}$ is  an integer sequence in  $ \mathcal{I}_k(\lambda)$ and  $\psi(\bsy{i})=w$. 
    Then Algorithm A
    transforms $w$ to $\bsy{i}$. 
\end{theorem}
\begin{proof}
    Let $w=b_1\cdots b_{n-k}\ a_1 \cdots a_k$. Since  $w \in \mathcal{R}_k(\lambda)$, we have that $b_1< b_2 < \cdots < b_{n-k}$, and in the process of applying the RSK algorithm to $w$, each term $a_i$ bumps a distinguished number 
    $t_i \in \{ b_1, \ldots, b_{n-k}, a_1, \ldots, a_{i-1}\}$  out of the first row. By the definition of row-insertion, 
    $t_i$ is the minimal integer that is larger than $a_i$ and not equal to $t_1, \ldots, t_{i-1}$. Hence, $t_i$ is 
    given by Step (2) of Algorithm A, and the insertion tableau of $\bsy{t}=(t_1, \ldots, t_k)$ under the RSK algorithm is $P(w)^*$, which is obtained from the 
     insertion tableau $P(w)$ of $w$ by removing the first row. 

    Let $P_0(\bsy{t})=\emptyset$ and $P_j(\bsy{t})$ be the insertion tableau of $(t_1,\ldots, t_j)$ under RSK for $j=1, \ldots, k$.    
    Let $(T^{(0)}, T^{(\frac12)}, T^{(1)}, \ldots, T^{(k)})$ be the sequence of  tableaux defined when applying $DI_n^k$ to  $\bsy{i}$.  
    By the definition of $\psi$, $P_j(\bsy{t})$ can be obtained from $T^{(j)}$ by removing the first row, and the first row of $T^{(j)}$ contains exactly integers in $[n] \setminus \mathrm{content}(P_j(\bsy{t}))$. 
    It follows that the shape of $T^{(j)}$ (resp. $T^{(j-\frac12)}$) is obtained from that of $P_j(\bsy{t})$ by adding a first row of $n-j$ (resp. $n-j-1$) boxes. 

    Recall that in the inverse map of $DI_n^k$, for $j=k, k-1, \ldots, 1$, we get the integer $i_j$ as the unique number $T^{(j)}= (i_j \xrightarrow[]{RSK} T^{(j-\frac12)})$.  From the definition of row-insertion, we see that $i_j$ is the largest  number in the first row of $T^{(j)}$ that is smaller than $t_j$, as described by Step (3) of Algorithm A.      
\end{proof}

\begin{corollary}
    Let $n \geq k+1$.  The set of integer sequences $\bsy{i} \in [n]^k$ with VT-index $k$ can be obtained by applying Algorithm A to permutations in $\mathcal{R}_k^n$. 
\end{corollary}

\begin{example}
    Let $n=4$ and $k=2$.  There are five integer sequences $(i_1, i_2) \in [4]^2$ with VT-index $2$, and $\mathcal{R}_2^4=\{ 1432, 2413, 2431, 3412, 3421\}$. By Algorithm A we get that the 
    set of integer sequences with VT-index $2$ are the union of $\mathcal{I}_2((2,2))=\{(2,2), (1,3)\}$ 
    and $\mathcal{I}_2((2,1,1))=\{ (3,2), (3,1), (2,1)\}$.     
\end{example}

Algorithm A can be easily inverted, giving a direct way to obtain $w=\psi(\bsy{i})$ for $\bsy{i} \in \mathcal{I}_k(\lambda)$. 
We now describe the inverse operations in Algorithm B. 

\begin{center}
\fbox{\begin{minipage}[t][1.1\height][c]{\dimexpr\textwidth-2\fboxsep-2\fboxrule\relax}
\noindent \textbf{Algorithm B: from  $\mathcal{I}_k(\lambda)$ to $\mathcal{R}_k(\lambda)$ } 
\begin{enumerate}
    \item Start with $\bsy{i}=(i_1, i_2, \ldots, i_k) \in \mathcal{I}_k(\lambda)$. 
    \item For $j=1$ to $k$, let 
    \[
    t_j := \min\{ x \in \mathbb{Z}: \ x > i_j \text{ and }  x \neq t_1, \ldots, t_{j-1}, i_{j+1}, \ldots, i_k\}.
    \]
    \item For $j=k$ to $1$, let 
    \[
     a_j:=\max\{ c \in \mathbb{Z}: \ c < t_j \text{ and } c \neq t_1, \ldots, t_{j-1}, a_{j+1}, \ldots, a_k\}. 
    \]
    \item Let $w=b_1 b_2  \cdots b_{n-k}\  a_1 a_2  \cdots a_k$, where $b_1 < b_2 < \cdots < b_{n-k}$ are the elements in 
    $[n]\setminus \{a_1, \ldots, a_k\}$. 
\end{enumerate}
\end{minipage}}
\end{center}

  Algorithm B can be visualized as operations in a rectangle grid. 
Given $\bsy{i}=(i_1, i_2, \ldots, i_k) \in \mathcal{I}_k(\lambda)$, for each $j=1, \ldots, k$ 
use diagonal shading 
\raisebox{-3pt}{\tikz{\fill[pattern=north east lines, pattern color=lightgray] (0,-1) rectangle (0.5,-1.5);\draw(0,-1) rectangle (0.5,-1.5)}} to shade 
all the cells 
$(i_j,j')$  
 for $1 \leq j' \leq j$. Starting from $j = 1$ to $j =k$, (i) find the lowest row $t_j$ above $(i_j,j)$ in column $j$ such that the cell $(t_j, j)$ is not shaded either diagonally \raisebox{-3pt}{\tikz{\fill[pattern=north east lines, pattern color=lightgray] (0,-1) rectangle (0.5,-1.5);\draw(0,-1) rectangle (0.5,-1.5)}} or horizontally \raisebox{-3pt}{\tikz{\fill[pattern=horizontal lines, pattern color=lightgray] (0,-1) rectangle (0.5,-1.5);\draw(0,-1) rectangle (0.5,-1.5)}}, and (ii) 
shade all the cells $(t_j,j')$ horizontally \raisebox{-3pt}{\tikz{\fill[pattern=horizontal lines, pattern color=lightgray] (0,-1) rectangle (0.5,-1.5);\draw(0,-1) rectangle (0.5,-1.5)}} for all $j \leq j'\leq k$.

Next,  remove  the diagonal shading but keep the horizontal ones. 
Starting from $ j= k$ to $j = 1$,  find the highest row $a_j \leq t_j$ 
such that $(a_j,j)$ is not shaded and row $a_j$ does not have a black cell. 
Color cell $(a_j,j)$ black.  Let $b_1 < b_2 < \cdots < b_{n-k} = n$ be the elements in $[n] \backslash \{a_1, a_2, \ldots, a_k\}$ in an increasing order. Then $w = b_1\cdots b_{n-k}\ a_1\cdots a_k \in \mathcal{R}_k(\lambda)$ and $\bsy{t}=(t_1, t_2, \ldots, t_k)$ is the corresponding bumping sequence. 

\begin{example}
Let $n =8$ and $k = 4$. 
We start with $\bsy{i} = (3,2,6,2) \in \mathcal{I}_4((4,2,2))$. The grids in the first row shows how to obtain the bumping sequence $\bsy{t} = (4,3,7,5)$, and the second row shows how to obtain 
$a_4=2$, $a_3=6$, $a_2=1$ and $a_1=3$,  which are represented by the black squares. From this we 
  obtain $w =45783162 \in \mathcal{R}_4((4,2,2))$.

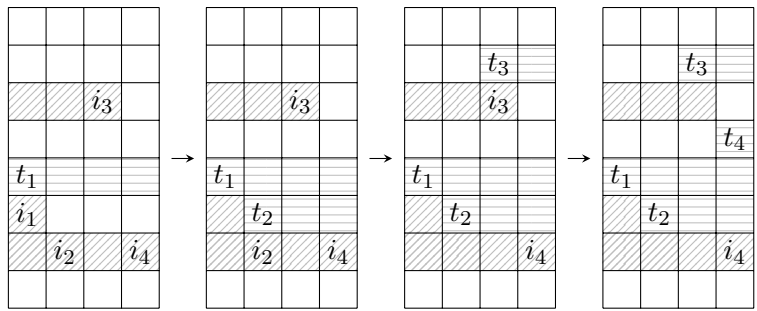
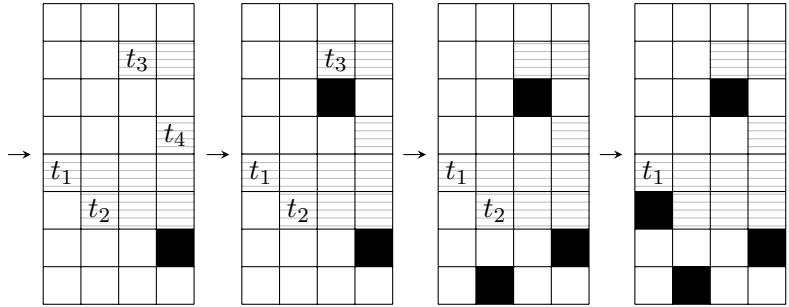
\input{image1_paper_2}

\end{example}

\medskip

Given an integer sequence $\bsy{i} \in [n]^k$, applying Algorithm B, if either  in Step (2), $t_j >n$ for some~$j$,  
 or in Step (3) $a_j \leq 0$ for some $j$, then  
 the VT-index of $\bsy{i}$ is not $k$. 
\begin{example} \label{ex: n4k2}
    Let $n=4$ and $k=2$. 
 \begin{enumerate}[(i)]
     \item For $\bsy{i}=(4,4)$,  Step (2) of Algorithm B gives $\bsy{t}=(5, 6) \not \in [n]^k$. Hence the VT-index of $(4,4)$ is not $2$. 
     \item For $\bsy{i} = (1,1)$, Step (2) gives $\bsy{t}=(2,3)$ and then Step (3) gives $a_2=1, a_1=0$. Hence the VT-index of $(1,1)$ is not $2$. 
\end{enumerate}
\end{example} 
However, there exists an integer sequence $\bsy{i} \in [n]^k$ that can pass through Algorithm B and yield  $\bsy{t}\in [n]^k$ and $w \in \mathcal{R}^n_k$, yet $\bsy{i}$ does not have VT-index $k$. For example, applying Algorithm B to $(1,2)$ we obtain 
     $\bsy{t}=(3,4)$ and $w=3412$, yet $\mathrm{vt}_4((1,2))=1$. In fact, we have  $\psi^{-1}(3412)=(2,2)$.

 On the other hand,  the correspondence between $w \in \mathcal{R}^n_k$ and the bumping sequence $\bsy{t}$ is one-to-one. 
More precisely, for a permutation $w=b_1 \cdots b_{n-k} \ a_1 \cdots a_k$ 
with $b_1< b_2 < \cdots < b_{n-k}$, $w \in \mathcal{R}^n_k$ if and only of Algorithm A(2) gives a sequence $\bsy{t}=(t_1, \ldots, t_k) \in [n]^k$. 
     This follows from the RSK algorithm and  Schensted's Theorem. 
     Conversely,  for a sequence $\bsy{t}=(t_1, \ldots, t_k) \in [n]^k$, if Algorithm B(3) gives a sequence $(a_1,  \ldots, a_k) \in [n]^k$, then $w=b_1 \cdots b_{n-k} \ a_1\cdots a_k \in \mathcal{R}^n_k$, where $b_1 < b_2 < \cdots < b_{n-k} $
     are the elements in $[n]\backslash \{a_1, \ldots, a_k\}$.   This follows from the fact that row-insertion is invertible.

The above argument implies that we can combine Algorithms A and B together to build a test to check whether 
a sequence $\bsy{i} \in [n]^k$ has 
 VT-index  $k$:
\begin{theorem}  \label{thm:test}
The VT-index of a sequence $\bsy{i} \in [n]^k$ is $k$ if and only if the following two conditions hold. 
\begin{enumerate}[(a)]
    \item  The sequences $\bsy{t}$ and $(a_1,a_2, \ldots, a_k)$   obtained by applying Algorithm B to $\bsy{i} $ are in $[n]^k$. 
    \item Applying Step (3) of Algorithm A  to $\bsy{t}$, we recover the  sequence $\bsy{i}$. 
\end{enumerate}
\end{theorem}

 Note that condition (b) of Theorem~\ref{thm:test} is not sufficient. 
 For instance, in (ii) of Example~\ref{ex: n4k2}, when $n=4$, $k=2$, and $\bsy{i}=(1,1)$, applying Step (2) of Algorithm B, we get $\bsy{t}=(2,3)$. Then applying Step (3) of Algorithm A to $(2,3)$, we recover $(1,1)$. Yet, $(1,1)$ has VT-index $1$. The missing step is Step (3) of Algorithm B.

In our algorithms  
 the three sequences, $\bsy{i}$,  $\bsy{t}$, and 
 $(a_1, \dots, a_k)$, are closely related. We know that the sequence $(a_1, \ldots, a_k)$ consists of the last $k$ terms of permutations in $\mathcal{R}_k^n$.  
In the next section we give an explicit characterization of the bumping sequence $\bsy{t}$.

\vspace{.5cm} 

\section{The bumping sequences and a reparking problem}\label{sec:bumping}

From the construction of Algorithms A and B in Section~\ref{section: max VT-index}, $\bsy{t}$ is a bumping sequence of 
    a permutation $w \in \mathcal{R}^n_k$ if and only if  $\bsy{t}$  yields  $(a_1,a_2,  \ldots, a_k) \in [n]^k$ in Algorithm B(3).
It raises the question:  Given a sequence $\bsy{t} \in [n]^k$ with distinct terms, under what condition does Algorithm B(3) yield a sequence $(a_1,a_2,  \ldots, a_k) \in [n]^k$?  We answer this question in Theorem~\ref{thm:criterion}.

For a sequence of distinct  integers $\bsy{a}=(a_1, a_2,\ldots, a_n)$, let $g(a_i)$ be the maximal length among all of the  increasing subsequences of $\bsy{a}$ that end at $a_i$.  
The following proposition is well-known. For example, see \cite[Chapter 8.4]{Aigner07}.  

\begin{proposition} \label{prop:first-row}
  Let $P(\bsy{a})$ be the insertion tableau of $\bsy{a}$ under the RSK algorithm. Then the entry in the $j$-th box of the first row of $P(\bsy{a})$ is $a_i$, where $i=\max\{k: g(a_k)=j\}$.  In particular, when inserting $a_i$, it bumps the term $a_t$ where $t$ is the largest index less than $i$ satisfying $g(a_t)=j$.   
\end{proposition}

\begin{theorem} \label{thm:criterion}
Let $n \geq k+1$ and $(t_1,\ldots, t_k) \in [n]^k$ be a  sequence of distinct integers. Let $t_{\max}= \max\{t_1,\ldots,t_k\}$ and $t_{\min} =  \min\{t_1,\ldots,t_k\}$. Let $m_1 < \cdots < m_j$ denote the elements in the set $\{t_{\min}, t_{\min}+1,\ldots, t_{\max}\} \setminus \{t_1,\ldots,t_k\}$. 
Let $l=\mathrm{is}(t_1, \ldots t_k, m_1, \dots, m_j) $  
be the length of the longest increasing subsequence of $(t_1,\ldots, t_k, m_1, \ldots, m_j)$. Then 
$(t_1,\ldots, t_k)$ is the bumping sequence of a permutation $w \in \mathcal{R}^n_k \subseteq \mathfrak{S}_n$ if and only if $t_{\min} > l-j$. 
\end{theorem}
\begin{proof}

Suppose that  $\bsy{t} =(t_1,\ldots,t_k) \in [n]^k$ is the bumping sequence corresponding to a permutation $w \in \mathcal{R}^n_k.$ By the definition of RSK, the sequence $\bsy{t}$ can be generated by applying Algorithm A(2) to $w$. 
We will construct two sequences, $\boldsymbol{u}$ and $\boldsymbol{t_m}$, and prove $t_{\min} > l-j$ by comparing the RSK insertion tableaux of these two sequences.
Consider the sequence 
\begin{equation} \label{def:u} 
\bsy{u} = (t_1, \ldots, t_k, 1,2,\ldots,t_{\min}-1,m_1,\ldots, m_j,t_{\max}+1,\ldots, n).
\end{equation}
Since $\bsy{t}$ is Knuth equivalent to the reading word of its RSK insertion tableau $P(\bsy{t})$, and 
the first row of the tableau $P(w)$ is exactly $(1,2,\ldots,t_{\min}-1,m_1,\ldots, m_j,t_{\max}+1,\ldots, n)$, we conclude that 
 $\bsy{u}$ shares the same RSK insertion tableau as the permutation  $w$. Moreover, the length of the first row of this insertion tableau is at least as long as the length of the first row of the insertion tableau of the sequence $\bsy{t_m} = (t_1, \ldots t_k,m_1,\ldots, m_j, t_{\max}+1,\ldots, n)$, which is $l+n-t_{\max}$ by Schensted's Theorem. Therefore  $n-k \geq l+ n-t_{\max}$. Noting that  $t_{\max}= t_{\min} +k+j-1$, we get 
 $t_{\min} \geq l-j+1 > l-j$.

Next, we show that if $t_{\min} > l-j$, there exists a sequence $w \in \mathcal{R}^n_k$ such that its bumping sequence is $\bsy{t}$. To accomplish this, we proceed with the following construction. Assume that $t_{\min} > l-j$. Consider the sequence $\bsy{u}$ as defined in Formula \eqref{def:u}.
Then the maximal size of increasing subsequences of $\bsy{u}$ is $\max\{t_{\min}-1+j+ n-t_{\max}, l+n-t_{\max}\} = t_{\min}-1+j+n-t_{\max}.$
By Proposition~\ref{prop:first-row},
the entries in the first row of the insertion tableau $P(\bsy{u})$ of $\bsy{u}$ are $1, 2, \ldots, t_{\min}-1, m_1,\ldots, m_j, t_{\max} +1, \ldots, n$. 
Consequently, the insertion tableau for the bumping sequence $\bsy{t}$ is positioned below the first row of $P(\bsy{u})$. 
We can then apply the inverse RSK algorithm to construct a sequence $w$ with the same insertion tableau as $P(\bsy{u})$.  
The recording tableau for $w$ has entries $1, 2, \ldots, n - k$ in the first row. Below the first row, the recording tableau is formed by adding $n-k$ to each entry of the recording tableau for the bumping sequence $(t_1, \ldots, t_k)$. Clearly, such a sequence $w$ belongs to $\mathcal{R}^n_k$.
\end{proof}

Comparing Algorithm A(2) and Algorithm B(3), we observe that there is a symmetry between the sequence $\bsy{a}=(a_1, \dots, a_k)$ and  $\bsy{t}$. Namely, $\bsy{t}$ is the bumping sequence of $\bsy{a}$ if and only if 
$\bsy{a}^{rc}$ is the bumping sequence of $\bsy{t}^{rc}$, where for a sequence $\bsy{x}=(x_1, \dots, x_k) \in [n]^k$,   \[\bsy{x}^{rc}:=(n+1-x_k, n+1-x_{k-1}, \ldots, n+1-x_1)\] is the reverse complement of $\bsy{x}$. 
Hence, Theorem~\ref{thm:criterion} gives the following characterization of permutations in $\mathcal{R}^n_k$. 

\begin{corollary} \label{cor:criterion2}
Let $n \geq k+1$ and $(a_1,a_2,\ldots, a_k) \in [n]^k$ be a  sequence of distinct integers. Let $a_{\max}= \max\{a_1,a_2,\ldots,a_k\}$ and $a_{\min} =  \min\{a_1,\ldots,a_k\}$. Let $n_1 < \cdots < n_j$ denote the elements in the set $\{a_{\min}, a_{\min}+1,\ldots, a_{\max}\} \setminus \{a_1,\ldots,a_k\}$. Let $l'=\mathrm{is}(n_1, \ldots, n_j, a_1,\ldots, a_k)$. 
Then 
$(a_1,\ldots, a_k)$ is the last $k$ entries  of a permutation $w \in \mathcal{R}^n_k \subseteq \mathfrak{S}_n$ if and only if $n+1-a_{\max} > l'-j$. 
\end{corollary}
\begin{proof}
    The sequence $\bsy{a}=(a_1, \ldots, a_k)$ is obtained from taking the last $k$ entries of a permutation $w \in \mathcal{R}^n_k $  if and only if $\bsy{a}^{rc}$ is a bumping sequence as described in Theorem~\ref{thm:criterion}. 
    Hence, $\min\{n+1-a_i\} > l'-j$, where $l'= 
\mathrm{is}(n+1-a_k, \ldots, n+1-a_1, n+1-n_j, \ldots, n+1-n_1)$, which is equal to  $\mathrm{is}(n_1, \ldots, n_j, a_1, \dots, a_k)$.    
\end{proof}

Theorem~\ref{thm:criterion} and Corollary~\ref{cor:criterion2} can be described  in terms of the following \emph{reparking problem}. 

Suppose that there are $k$ cars parked  on a street with $n$
parking spots numbered $1$ through $n$ from 
the left to right, where car $i$ is parked at position $x_i$. 
Assume $n_1< \cdots  < n_j$ are the empty spots between $x_{\min}$ and $x_{\max}$. 
\begin{enumerate}
    \item Starting from car 1, for $j=1, 2, \ldots, k$, we ask car $j$ one-by-one to move to the closest available spot on its right.  Reparking is successful if all of the cars can find a new spot on the street. 

     Then the reparking is successful if and only if the number of empty spots at the end of the street is at least   $l'-j$, where $l'=\mathrm{is}(n_1, \ldots, n_j, x_1, \ldots, x_k)$.

    \item Starting from car $k$, for $j=k, k-1, \ldots, 1$, we
ask  car $j$ one-by-one to move to the closest available spot on its left. 
Reparking is successful if  all of the cars can find a new spot on the street.

Then the reparking is successful  if and only if the number of empty spots at the beginning of the street is 
at least $l-j$, where $l=\mathrm{is}( x_1, \ldots, x_k, n_1, \ldots, n_j  )$.  
\end{enumerate}

\begin{example}
    Let $k=5$ and assume that cars are parked at positions $\bsy{x}=(3, 2, 5, 8, 9)$. Then $j=3$ and 
    $(n_1, n_2, n_3)=(4, 6,7)$. 
\begin{enumerate}
    \item For $j=1, 2, \ldots, k$, move car $j$ one-by-one to the closest available spot on its right.
    Then the cars are moved to spots $(4, 3, 6, 10, 11)$. The cars occupy two new spots on the right. Note that 
    $l'=\mathrm{is}(4, 6, 7, 3, 2, 5, 8, 9)=5$ and $2=l'-j$. 
    \item For $j=k, k-1, \ldots, 1$, we move car $j$ one-by-one to  the closest available spot on its left.
    Then the cars are moved to spots $(2, 1, 4, 6, 7)$, where car 1 moves to spot 1, car 2 moves to spot 2, car 3 moves to spot 4, car 4 moves to spot 6, and car 5 moves to spot 7. 
        The cars occupy one new spots on the left.  Note that 
$l=\mathrm{is}(3,2,5,8,9,4,6,7)= 4$ and $l-j=1$. 
\end{enumerate}
    
\end{example}

\section{Final Remark} \label{sec:final} 

Recently, there is another bijective proof of Identity \eqref{eq: Identity1} constructed by Colmenarejo, Orellana, Saliola, Schilling, and Zabrocki~\cite{COSSZ}, which we call the COSSZ bijection. This bijection 
maps integer sequences in $[n]^k$ to pairs of tableaux of the same  shape $\lambda \vdash n$
with one being a standard Young tableau of content $[n]$ and the other being a  standard multiset tableaux of content $[k]$. 
A standard multiset tableau of shape $\lambda$ is a filling of the Young diagram of $\lambda$ with disjoint subsets of $[k]$ that is increasing along each row and each column, where disjoint subsets are ordered by 
their maximum elements, i.e., for two disjoint subsets $A$ and $B$ of $[k]$, $A \prec B$ if and only if $\max(A) \leq \max(B)$. 
In particular, $\emptyset \prec A$ if $A$ 
is non-empty.

  Under the COSSZ bijection,   each integer sequence $\bsy{i} \in [n]^k$ is mapped to pairs of tableaux of the same shape $\lambda \vdash n$. We can again ask: for which $\bsy{i}$ does the shape $\lambda$ satisfy $\lambda_1=n-k$?   We show that such integer sequences  can again be characterized by the set 
 $\mathcal{R}^n_k \subseteq \mathfrak{S}_n$, in a more direct way compared to the integer sequences under $DI_n^k$.  Note that the COSSZ bijection  follows the same idea as in \cite{CDDSY07}, using  
 only the RSK insertion (Knuth's version) but not jeu de taquin.  The difference between the result here and those in Section 4 shows the intrinsic complexity and the subtle relation between the row-insertion and jeu de taquin.

Below is the description of the COSSZ bijection, which is defined for all integers $n \geq 1$  and $k \geq 0$. 
Let  $\bsy{i}=(i_1, i_2, \ldots, i_k) \in [n]^k$ be an integer sequence. 
\begin{enumerate}[(a)]
    \item Let $M_r=\{ j: i_j=r\}$ be the positions of integer $r$ in the sequence $\bsy{i}$. 
       Then the non-empty $M_r$'s form a set partition $P$ of $[k]$. 
       Assume that $P$ has $t$ blocks.  

    \item Form a two-line  array 
    \[
    A= \begin{pmatrix} 
        a_1 & a_2 & \cdots & a_n \\ 
        b_1 & b_2 & \cdots & b_n 
    \end{pmatrix} 
    \]
         where $a_1 \leq  a_2 \leq  \cdots \leq a_n$ is the non-decreasing rearrangement of $\{ \emptyset^{n-t} \} \cup P$, ordered by their largest elements, and 
         $w= b_1b_2\cdots b_n$ is a permutation of length $n$ such that 
         \begin{enumerate}[(i)]
             \item $b_1 < b_2 < \cdots < b_{n-t}$ are the integers  $x \in [n]$ such that $M_x=\emptyset$. That is,  
              $M_{b_1}=  \cdots = M_{b_{n-t}}=\emptyset$. 
              
            \item   $a_i=M_{b_i}$  if $M_{b_i} \neq \emptyset$. 
              \end{enumerate}
          \item Apply the RSK algorithm to the two-line array $A$ to obtain a pair of tableaux
           $(S, T)$ of the same shape $\lambda$, where $S \in \mathcal{SYT}(\lambda)$ is the insertion tableau of $b_1b_2 \cdots b_n$, and $T$ is  the recording tableau, whose entries are subsets $M_r$'s, for 
           $r=1, 2, \ldots, n$.
   \end{enumerate}
The image of $\bsy{i}$ is then $(S, T)$. 

      \begin{example}
           Let $n=6$, $k=4$, and  $\bsy{i}=(3,2,6,2)$.  Then  $M_2=\{2, 4\}$, $M_3=\{1\}$, $M_6=\{3\}$ and $M_1=M_4=M_5=\emptyset$. Hence 
           the two-line array is 
           \[
           A=\begin{pmatrix}
               \emptyset & \emptyset  & \emptyset  & \{ 1\}  & \{ 3 \} & \{ 2,4\}  \\ 
               1 & 4 &  5 &   3  & 6 & 2 
           \end{pmatrix}.
           \]
          Applying Knuth's RSK, we obtain that the image of $\bsy{i}$ under the COSSZ bijection  is 
          \[
         S=\ytableaushort{1256,3,4}, \qquad T= \ytableaushort{\emptyset\emptyset \emptyset 3,1,{24}}.
          \]        
           \end{example}
     
     Assume $n \geq k+1$.  We consider when the image of $\bsy{i}$ has a shape $\lambda \vdash n$ with $\lambda_1=n-k$. 
In that case, each entry not in the first row of $T$ must be a non-empty subset of $[k]$, hence the partition formed by 
non-empty $M_r$'s has at least $k$ blocks.  But it is a set partition of $[k]$, hence 
        each non-empty $M_r$ must be a singleton block and   the entries of the first row of $T$ 
      must all be the empty set. 
  It follows that   $a_1=a_2=\cdots =a_{n-k}= \emptyset$ and $a_{n-k+j}=\{j\}$ for $j \in [k]$. Therefore 
  $b_1 < b_2 < \cdots < b_{n-k}$ and $b_{n-k+j}=i_j$, i.e.,  $\bsy{i}$ is exactly $(b_{n-k+1}, b_{n-k+2},  \ldots,  b_n)$, the last $k$ integers of the second row of the two-line array. 

   Note that $w= b_1 b_2 \cdots b_n$ is a permutation of length $n$. 
   The shape of $P$ and $T$ is $\lambda$ implies that $\mathrm{is}(b_1b_2 \cdots b_n) = \lambda_1= n-k$. 
    Hence,  $w \in \mathcal{R}^n_k$.  In summary, under the COSSZ bijection,  an integer sequence $\bsy{i} \in [n]^k$ is mapped to a pair  of tableaux of shape $\lambda \vdash n$ with $\lambda_1=n-k$ 
        if and only if  $\bsy{i}$ consists of the last $k$ entries of a permutation 
       $w \in \mathcal{R}^n_k$.

       \begin{example}
           Let $n=8$, $k=4$, and  $\bsy{i}=(3,1,6,2)$.  Then $M_1=\{2\}$, $M_2=\{4\}$, $M_3=\{1\}$, $M_6=\{3\}$ and $M_4=M_5=M_7=M_8=\emptyset$. Hence 
           the two-line array is 
           \[
           A=\begin{pmatrix}
               \emptyset & \emptyset & \emptyset & \emptyset & \{1\} & \{2 \} & \{ 3\}  & \{4\} \\ 
               4 & 5 & 7 & 8 & 3 & 1 & 6 & 2 
           \end{pmatrix}. 
           \]
          Applying Knuth's RSK, we obtain that the image of $\bsy{i}$ under the COSSZ bijection  is 
          \[
         S=\ytableaushort{1268,35,47}, \qquad T= \ytableaushort{\emptyset \emptyset \emptyset \emptyset,13,24}
          \]
           and  the shape of $S$ and $T$  is $\lambda=(4,2,2)$, which satisfies $\lambda_1=n-k$. 
          Observe that $\bsy{i}$ consists of the last four letters of the permutation $w = 45783162$, 
           which is in $\mathcal{R}^8_4$.

       \end{example}

\section*{Acknowledgements}
This material is based upon work supported by the National Science Foundation under Grant No. DMS-1929284 while the authors was in residence at the Institute for Computational and Experimental Research in Mathematics (ICERM) in Providence, RI, during the Discrete Optimization: Mathematics,  Algorithms, and Computation semester program.
We thank ICERM for facilitating and supporting this research project. We are also grateful to James Sundstrom for providing help on computer programming. 

The second author is supported in part by an EDGE Karen Uhlenbeck fellowship, the third author is supported by a PSC-CUNY award, and the fourth author is supported in part by the Simons Collaboration Grant for Mathematics 704276.

\end{document}

%% file: image_revised.tex
\[\begin{tikzpicture}

\fill[pattern=north east lines, pattern color=lightgray] (0,2.5) rectangle (1.5,3);
\fill[pattern=north east lines, pattern color=lightgray] (0,1) rectangle (0.5,1.5);
\fill[pattern=north east lines, pattern color=lightgray] (0,0) rectangle (1,0.5);
\fill[pattern=north east lines, pattern color=lightgray] (0,0.5) rectangle (2,1);

\draw (0.25,1.25) node{$a_1$};
\draw (0.75,0.25) node{$a_2$};
\draw (1.25,2.75) node{$a_3$};
\draw (1.75,0.75) node{$a_4$};

\fill[pattern=horizontal lines light gray] (0,1.5) rectangle (2,2);
\draw (0.25,1.75) node{$t_1$};
\draw[step=0.5cm,black, thin] (0,0) grid (2,4);
\end{tikzpicture} \hspace{.05in} \begin{tikzpicture}[anchor=base, baseline]
\draw [-stealth](0,2) -- (0.3,2);
\end{tikzpicture}\hspace{.05in} 
\begin{tikzpicture}

\fill[pattern=north east lines, pattern color=lightgray] (0,2.5) rectangle (1.5,3);
\fill[pattern=north east lines, pattern color=lightgray] (0,1) rectangle (0.5,1.5);
\fill[pattern=north east lines, pattern color=lightgray] (0,0) rectangle (1,0.5);
\fill[pattern=north east lines, pattern color=lightgray] (0,0.5) rectangle (2,1);

\draw (0.75,0.25) node{$a_2$};
\draw (1.25,2.75) node{$a_3$};
\draw (1.75,0.75) node{$a_4$};
\fill[pattern=horizontal lines light gray] (0,1.5) rectangle (2,2);
\draw (0.25,1.75) node{$t_1$};
\fill[pattern=horizontal lines light gray] (0.5,1) rectangle (2,1.5);
\draw (0.75,1.25) node{$t_2$};
\draw[step=0.5cm,black, thin] (0,0) grid (2,4);
\end{tikzpicture}\hspace{.05in}
\hspace{.05in} 
\begin{tikzpicture}[anchor=base, baseline]
\draw [-stealth](0,2) -- (0.3,2);
\end{tikzpicture}\hspace{.05in} \begin{tikzpicture}

\fill[pattern=north east lines, pattern color=lightgray] (0,2.5) rectangle (1.5,3);
\fill[pattern=north east lines, pattern color=lightgray] (0,1) rectangle (0.5,1.5);
\fill[pattern=north east lines, pattern color=lightgray] (0,0) rectangle (1,0.5);
\fill[pattern=north east lines, pattern color=lightgray] (0,0.5) rectangle (2,1);

\draw (1.25,2.75) node{$a_3$};
\draw (1.75,0.75) node{$a_4$};
\fill[pattern=horizontal lines light gray] (0,1.5) rectangle (2,2);
\draw (0.25,1.75) node{$t_1$};
\fill[pattern=horizontal lines light gray] (0.5,1) rectangle (2,1.5);
\draw (0.75,1.25) node{$t_2$};
\fill[pattern=horizontal lines light gray] (1,3) rectangle (2,3.5);
\draw (1.25,3.25) node{$t_3$};
\draw[step=0.5cm,black, thin] (0,0) grid (2,4);
\end{tikzpicture} \hspace{.05in}\begin{tikzpicture}[anchor=base, baseline]
\draw [-stealth](0,2) -- (0.3,2);
\end{tikzpicture}\hspace{.05in}  \begin{tikzpicture}

\fill[pattern=north east lines, pattern color=lightgray] (0,2.5) rectangle (1.5,3);
\fill[pattern=north east lines, pattern color=lightgray] (0,1) rectangle (0.5,1.5);
\fill[pattern=north east lines, pattern color=lightgray] (0,0) rectangle (1,0.5);
\fill[pattern=north east lines, pattern color=lightgray] (0,0.5) rectangle (2,1);

\draw (1.75,0.75) node{$a_4$};
\fill[pattern=horizontal lines light gray] (0,1.5) rectangle (2,2);
\draw (0.25,1.75) node{$t_1$};
\fill[pattern=horizontal lines light gray] (0.5,1) rectangle (2,1.5);
\draw (0.75,1.25) node{$t_2$};
\fill[pattern=horizontal lines light gray] (1,3) rectangle (2,3.5);
\draw (1.25,3.25) node{$t_3$};
\fill[pattern=horizontal lines light gray] (1.5,2) rectangle (2,2.5);
\draw (1.75,2.25) node{$t_4$};
\draw[step=0.5cm,black, thin] (0,0) grid (2,4);
\end{tikzpicture} \hspace{.05in} \begin{tikzpicture}[anchor=base, baseline]
\draw [-stealth](0,2) -- (0.3,2);
\end{tikzpicture}\hspace{.05in}  \begin{tikzpicture}

\fill[pattern=horizontal lines light gray] (0,1.5) rectangle (2,2);
\draw (0.25,1.75) node{$t_1$};
\fill[pattern=horizontal lines light gray] (0.5,1) rectangle (2,1.5);
\draw (0.75,1.25) node{$t_2$};
\fill[pattern=horizontal lines light gray] (1,3) rectangle (2,3.5);
\draw (1.25,3.25) node{$t_3$};
\fill[pattern=horizontal lines light gray] (1.5,2) rectangle (2,2.5);
\draw (1.75,2.25) node{$t_4$};
\fill[black] (0,1) rectangle (0.5,1.5);
\fill[black] (0.5,0.5) rectangle (1,1);
\fill[black] (1,2.5) rectangle (1.5,3);
\fill[black] (1.5,0.5) rectangle (2,1);
\draw[step=0.5cm,black, thin] (0,0) grid (2,4);
\end{tikzpicture}.\]

%% file: keyforfigures.tex
\noindent\fbox{\begin{minipage}[t][1\height][l]{\dimexpr\textwidth-2\fboxsep-2\fboxrule\relax}
\[\begin{tikzpicture}
\draw (-2,.25) node{Key:};
\fill[black]  (0,0) rectangle (0.5,.5);
 \fill[pattern=north east lines, pattern color=lightgray] (3,0) rectangle (3.5,.5);
 \draw[black, thin] (3.0,0)--(3.5,0)--(3.5,.5)--(3,.5)--(3,0);
 \fill[pattern=horizontal lines light gray] (7,0) rectangle (7.5,.5);
 \draw[black, thin] (7.0,0)--(7.5,0)--(7.5,.5)--(7,.5)--(7,0);
\draw (1.3,.25) node{= black};
\draw (4.6,.25) node{= diagonal};
\draw (8.7,.25) node{= horizontal};
\end{tikzpicture}
\]
\end{minipage}
}

%% file: image1_paper_2.tex
\[\begin{tikzpicture}
\fill[thin, pattern=horizontal lines, pattern color=lightgray] (0,1.5) rectangle (2,2);
\draw (0.25,1.75) node{$t_1$};
\fill[pattern=north east lines, pattern color=lightgray] (0,1) rectangle (0.5,1.5);
\fill[pattern=north east lines, pattern color=lightgray] (0,0.5) rectangle (1,1);
\fill[pattern=north east lines, pattern color=lightgray] (0,2.5) rectangle (1.5,3);
\fill[pattern=north east lines, pattern color=lightgray] (0,0.5) rectangle (2,1);
\draw (0.25,1.25) node{$i_1$};
\draw (0.75,0.75) node{$i_2$};
\draw (1.25,2.75) node{$i_3$};
\draw (1.75,0.75) node{$i_4$};

\draw[step=0.5cm,black, thin] (0,0) grid (2,4);
\end{tikzpicture}\hspace{.05in}\begin{tikzpicture}[anchor=base, baseline]
\draw [-stealth](0,2) -- (0.3,2);
\end{tikzpicture}\hspace{.05in}
\begin{tikzpicture}
\fill[thin, pattern=horizontal lines, pattern color=lightgray] (0,1.5) rectangle (2,2);
\fill[thin, pattern=horizontal lines, pattern color=lightgray] (0.5,1) rectangle (2,1.5);
\draw (0.25,1.75) node{$t_1$};
\draw (0.75,1.25) node{$t_2$};
\fill[pattern=north east lines, pattern color=lightgray] (0,1) rectangle (0.5,1.5);
\fill[pattern=north east lines, pattern color=lightgray] (0,0.5) rectangle (1,1);
\fill[pattern=north east lines, pattern color=lightgray] (0,2.5) rectangle (1.5,3);
\fill[pattern=north east lines, pattern color=lightgray] (0,0.5) rectangle (2,1);
\draw (0.75,0.75) node{$i_2$};
\draw (1.25,2.75) node{$i_3$};
\draw (1.75,0.75) node{$i_4$};

\draw[step=0.5cm,black, thin] (0,0) grid (2,4);
\end{tikzpicture}\hspace{.05in}\begin{tikzpicture}[anchor=base, baseline]
\draw [-stealth](0,2) -- (0.3,2);
\end{tikzpicture}\hspace{.05in}
\begin{tikzpicture}
\fill[thin, pattern=horizontal lines, pattern color=lightgray] (0,1.5) rectangle (2,2);
\fill[thin, pattern=horizontal lines, pattern color=lightgray] (0.5,1) rectangle (2,1.5);
\fill[thin, pattern=horizontal lines, pattern color=lightgray] (1,3) rectangle (2,3.5);
\draw (0.25,1.75) node{$t_1$};
\draw (0.75,1.25) node{$t_2$};
\draw (1.25,3.25) node{$t_3$};
\fill[pattern=north east lines, pattern color=lightgray] (0,1) rectangle (0.5,1.5);
\fill[pattern=north east lines, pattern color=lightgray] (0,0.5) rectangle (1,1);
\fill[pattern=north east lines, pattern color=lightgray] (0,2.5) rectangle (1.5,3);
\fill[pattern=north east lines, pattern color=lightgray] (0,0.5) rectangle (2,1);
\draw (1.25,2.75) node{$i_3$};
\draw (1.75,0.75) node{$i_4$};

\draw[step=0.5cm,black, thin] (0,0) grid (2,4);
\end{tikzpicture}\hspace{.05in}\begin{tikzpicture}[anchor=base, baseline]
\draw [-stealth](0,2) -- (0.3,2);
\end{tikzpicture}\hspace{.05in}
\begin{tikzpicture}
\fill[thin, pattern=horizontal lines, pattern color=lightgray] (0,1.5) rectangle (2,2);
\fill[thin, pattern=horizontal lines, pattern color=lightgray] (0.5,1) rectangle (2,1.5);
\fill[thin, pattern=horizontal lines, pattern color=lightgray] (1,3) rectangle (2,3.5);
\fill[thin, pattern=horizontal lines, pattern color=lightgray] (1.5,2) rectangle (2,2.5);
\draw (0.25,1.75) node{$t_1$};
\draw (0.75,1.25) node{$t_2$};
\draw (1.25,3.25) node{$t_3$};
\draw (1.75,2.25) node{$t_4$};
\fill[pattern=north east lines, pattern color=lightgray] (0,1) rectangle (0.5,1.5);
\fill[pattern=north east lines, pattern color=lightgray] (0,0.5) rectangle (1,1);
\fill[pattern=north east lines, pattern color=lightgray] (0,2.5) rectangle (1.5,3);
\fill[pattern=north east lines, pattern color=lightgray] (0,0.5) rectangle (2,1);
\draw (1.75,0.75) node{$i_4$};

\draw[step=0.5cm,black, thin] (0,0) grid (2,4);
\end{tikzpicture}\]
\[\begin{tikzpicture}[anchor=base, baseline]
\draw [-stealth](0,2) -- (0.3,2);
\end{tikzpicture}\hspace{.05in}\begin{tikzpicture}
\fill[thin, pattern=horizontal lines, pattern color=lightgray] (0,1.5) rectangle (2,2);
\fill[thin, pattern=horizontal lines, pattern color=lightgray] (0.5,1) rectangle (2,1.5);
\fill[thin, pattern=horizontal lines, pattern color=lightgray] (1,3) rectangle (2,3.5);
\fill[thin, pattern=horizontal lines, pattern color=lightgray] (1.5,2) rectangle (2,2.5);
\draw (0.25,1.75) node{$t_1$};
\draw (0.75,1.25) node{$t_2$};
\draw (1.25,3.25) node{$t_3$};
\draw (1.75,2.25) node{$t_4$};
\fill[black] (1.5,0.5) rectangle (2,1);
\draw[step=0.5cm,black, thin] (0,0) grid (2,4);
\end{tikzpicture}\hspace{.05in}\begin{tikzpicture}[anchor=base, baseline]
\draw [-stealth](0,2) -- (0.3,2);
\end{tikzpicture}\hspace{.05in}
\begin{tikzpicture}
\fill[thin, pattern=horizontal lines, pattern color=lightgray] (0,1.5) rectangle (2,2);
\fill[thin, pattern=horizontal lines, pattern color=lightgray] (0.5,1) rectangle (2,1.5);
\fill[thin, pattern=horizontal lines, pattern color=lightgray] (1,3) rectangle (2,3.5);
\fill[thin, pattern=horizontal lines, pattern color=lightgray] (1.5,2) rectangle (2,2.5);
\draw (0.25,1.75) node{$t_1$};
\draw (0.75,1.25) node{$t_2$};
\draw (1.25,3.25) node{$t_3$};

\fill[black] (1,2.5) rectangle (1.5,3);
\fill[black] (1.5,0.5) rectangle (2,1);
\draw[step=0.5cm,black, thin] (0,0) grid (2,4);
\end{tikzpicture}\hspace{.05in}\begin{tikzpicture}[anchor=base, baseline]
\draw [-stealth](0,2) -- (0.3,2);
\end{tikzpicture}\hspace{.05in}
\begin{tikzpicture}
\fill[thin, pattern=horizontal lines, pattern color=lightgray] (0,1.5) rectangle (2,2);
\fill[thin, pattern=horizontal lines, pattern color=lightgray] (0.5,1) rectangle (2,1.5);
\fill[thin, pattern=horizontal lines, pattern color=lightgray] (1,3) rectangle (2,3.5);
\fill[thin, pattern=horizontal lines, pattern color=lightgray] (1.5,2) rectangle (2,2.5);
\draw (0.25,1.75) node{$t_1$};
\draw (0.75,1.25) node{$t_2$};

\fill[black] (1,2.5) rectangle (1.5,3);
\fill[black] (1.5,0.5) rectangle (2,1);
\fill[black] (0.5,0) rectangle (1,0.5);
\draw[step=0.5cm,black, thin] (0,0) grid (2,4);
\end{tikzpicture}\hspace{.05in}\begin{tikzpicture}[anchor=base, baseline]
\draw [-stealth](0,2) -- (0.3,2);
\end{tikzpicture}\hspace{.05in}
\begin{tikzpicture}
\fill[thin, pattern=horizontal lines, pattern color=lightgray] (0,1.5) rectangle (2,2);
\fill[thin, pattern=horizontal lines, pattern color=lightgray] (0.5,1) rectangle (2,1.5);
\fill[thin, pattern=horizontal lines, pattern color=lightgray] (1,3) rectangle (2,3.5);
\fill[thin, pattern=horizontal lines, pattern color=lightgray] (1.5,2) rectangle (2,2.5);
\draw (0.25,1.75) node{$t_1$};
\fill[black] (0,1) rectangle (0.5,1.5);
\fill[black] (1,2.5) rectangle (1.5,3);
\fill[black] (1.5,0.5) rectangle (2,1);
\fill[black] (0.5,0) rectangle (1,0.5);
\draw[step=0.5cm,black, thin] (0,0) grid (2,4);
\end{tikzpicture}\]

%% file: On_the_Correspondence_Between_Integer_Sequences_and_Vacillating_Tableaux.bbl
\begin{thebibliography}{10}

\bibitem{Aigner07} M. Aigner. A Course in Enumeration. 
Graduate Text in Mathematics, Volume 238. Springer, 2007. 

\bibitem{BHH17} 
G.~Benkart, T.~Halverson, and N.~Harman, Dimensions of irreducible modules for partition algebras and tensor power multiplicities for symmetric and alternating groups, 
\emph{J. Algebr. Combin}. {\bf 46}(1) (2017), 77-108. 





\bibitem{BHPYZ1} Z.~Berikkyzy,  P.~E.~Harris, A.~Pun, C.~Yan, and C.~Zhao. 
On the limiting vacillating tableaux for integer sequences. 
\emph{J.  Comb.}, to appear.  

\bibitem{BHPYZ2} Z.~Berikkyzy,  P.~E.~Harris, A.~Pun, C.~Yan, and C.~Zhao. 
Combinatorial identities for vacillating tableaux. 
\emph{Integers}, to appear. 


\bibitem{CDDSY07} 
W.Y.C.~Chen, E.Y.P.~Deng, R.R.X.~Du, R.P.~Stanley and C.H.~Yan.
{Crossings and nestings of matchings and set partitions}.  
\emph{Trans. Amer. Math. Soc}. {\bf 359}(4)  (2007), 1555--1575.

 
\bibitem{COSSZ} L.~Colmenarejo, R.~Orellana,  F.~Saliola, A.~Schilling, and M.~Zabrocki. An insertion algorithm of multiset partitions with applications to diagram algebras. 
\emph{J. Algebra}, {\bf 557}(1) (2020),  97-128. 

\bibitem{DuYan24} 
R.~Du and C.H.~Yan.  On the enumeration of vacillating tableaux (in Chinese). \emph{Sci Sin Math}, 2025,
55: 1–11. 

\bibitem{GP2020} T.~Guo and S.~Poznanovi\'{c}. Hecke insertion and maximal increasing and decreasing sequences in fillings of stack polyominoes. 
\emph{J. Combin. Theory Ser A}, {\bf 176} (2020), 105304.



\bibitem{HL05} T.~Halverson and  T.~Lewandowski. 
RSK insertion for set partitions and diagram algebras.
\emph{Electron J. Combin.} {\bf 11} (2004/2006), \#R24, 24pp. 




\bibitem{Knuth70} D.E.~Knuth. Permutations, matrices, and generalized Young tableaux. 
\emph{Pacific J.  Math.}  {\bf 34}(3) (1970), 709-- 727. 

\bibitem{Kratten06} C.~Krattenthaler. Growth diagrams, and increasing and decreasing chains in fillings of Ferrers shapes. 
\emph{Adv. in Appl. Math.}, {\bf 37 }(2006), 404–431.

\bibitem{Handbook10} C.~Krattenthaler. Lattice path enumeration. 
  Chapter 10 in Handbook of Enumerative Combinatorics, edited by Miklos Bona. Chapman and Hall/CRC, 2015. ISBN 9781482220858.


\bibitem{Kratten23} C.~Krattenthaler. 
Identities for vacillating tableaux via growth diagrams, preprint. 
\href{https://arxiv.org/abs/2304.07657}{arXiv:2304.07657}. 

\bibitem{Rubey11} M.~Rubey. Increasing and decreasing sequences in fillings of moon polyominoes. 
\emph{Adv. in Appl. Math.} {\bf 47} (2011), 57–87.


\bibitem{MR98} P.P.~Martin and G.~Rollet, The Potts model representation and a Robinson–Schensted correspondence
for the partition algebra, \emph{Compos. Math}. {\bf 112} (1998), 237–254.





\bibitem{Robinson38} 
G. de B. Robinson, On the representations of the symmetric group, 
\emph{Amer. J. Math}. {\bf 60}(3) (1938),  745–760. 


 
\bibitem{Sagan01} B.~Sagan. the Symmetric Groups, Representations, Combinatorial Algorithms, and Symmetric Functions, second edition, Springer-Verlag, New York, 2001. 


\bibitem{EC2} R.P.~Stanley. Enumerative Combinatorics, Vol. 2, Cambridge University Press, Cambridge, 1999. 



\bibitem{Sch61} C.~Schensted, Longest increasing and decreasing subsequences. \emph{Canad. J. Math.} {\bf 13} (1961), 179--191. 


\end{thebibliography}
